\documentclass[11pt,a4paper]{amsart}
\usepackage{amsmath, graphicx, subfigure, epsfig,amssymb,cite}
\usepackage{xcolor,mathtools}
\usepackage[normalem]{ulem}
\numberwithin{equation}{section}

\usepackage{hyperref}
\hypersetup{
    colorlinks=true,
    linkcolor=blue,
    filecolor=magenta,      
    urlcolor=cyan,
}

\newcommand{\ssp}{\text{sing\,supp}}

\newcommand{\im}{\text{Im}}
\newcommand{\re}{\text{Re}}
\newcommand{\e}{\epsilon}
\newcommand{\ga}{\gamma}
\newcommand{\de}{\delta}
\newcommand{\br}{\mathbb{R}}

\newcommand{\BN}{\mathbb{N}}

\newcommand{\pa}{\partial}
\newcommand{\bt}{\beta}
\newcommand{\al}{\alpha}
\newcommand{\la}{\lambda}

\newcommand{\coi}{C_0^{\infty}}
\newcommand{\ioi}{\int_0^{\infty}}

\newcommand{\be}{\begin{equation}}
\newcommand{\ee}{\end{equation}}
\def\bs#1\es{
    \begin{equation}\begin{split}
    #1
    \end{split}\end{equation}
}
\newcommand{\dd}{\text{d}}
\newcommand{\fp}{\mathfrak p}
\newcommand{\fr}{f^{\text{rec}}}

\newcommand{\chl}{\check\la}
\newcommand{\chm}{\check\mu}

\newcommand{\bma}{\begin{pmatrix}}
\newcommand{\ema}{\end{pmatrix}}

\newcommand{\CL}{\mathcal{L}}

\newcommand{\CU}{\mathcal U}

\newcommand{\s}{\mathcal S}

\newcommand{\CF}{\mathcal F}
\newcommand{\CH}{\mathcal H}

\newcommand{\CD}{\mathcal D}

\newtheorem{theorem}{Theorem}[section]

\theoremstyle{definition}

\newtheorem{definition}[theorem]{Definition}
\newtheorem{remark}[theorem]{Remark}

\begin{document}

\title[Beam hardening streaks]{Analysis of beam hardening streaks in tomography}
\author[A Katsevich]{Alexander Katsevich}

\begin{abstract} 
The mathematical foundation of X-ray CT is based on the assumption that by measuring the attenuation of X-rays passing through an object, one can recover the integrals of the attenuation coefficient $\mu(x)$ along a sufficiently rich family of lines $L$, $\int_L \mu(x) \text{d} x$. This assumption is inaccurate because the energy spectrum of an X-ray beam in a typical CT scanner is wide. At the same time, the X-ray attenuation coefficient of most materials is energy-dependent, and this dependence varies among materials. Thus, reconstruction from X-ray CT data is a nonlinear problem. If the nonlinear nature of CT data is ignored and a conventional linear reconstruction formula is used, which is frequently the case, the resulting image contains beam-hardening artifacts such as streaks. In this work, we describe the nonlinearity of CT data using the conventional model accepted by all CT practitioners. Our main result is the characterization of streak artifacts caused by nonlinearity. We also obtain an explicit expression for the leading singular behavior of the artifacts. Finally, a numerical experiment is conducted to validate the theoretical results.
\end{abstract}

\maketitle

\section{Introduction}\label{sec:intro}

X-ray Computer Tomography (CT) is a ubiquitous imaging tool used in numerous applications, including medical imaging, nondestructive testing, industrial metrology, and security scanning. Let $\mu(x)$ denote the X-ray attenuation coefficient at point $x$ inside an object. The mathematical foundation of CT is based on the assumption that by measuring the attenuation of X-rays passing through the object, one can recover the integrals of $\mu(x)$ along a sufficiently rich family of lines $L$, $\int_L \mu(x)\dd x$. 

This assumption is usually inaccurate because the vast majority of conventional X-ray sources generate X-ray beams with a wide spectrum supported on an interval $[0,E_{\text{mx}}]$ \cite[section 4.2]{Podgorsak2010}, \cite{Forghani2017a}. Here $E_{\text{mx}}$ is the maximum energy of photons in the beam. At the same time, the X-ray attenuation coefficient of the vast majority of materials depends on energy, and this dependence differs for different materials. Hence, it is more appropriate to describe the attenuation coefficient as a function of position $x$ and energy $E$, that is, $\mu(x,E)$. Thus, using Beer's law, the CT data, $g(L)$, can be written in the form
\bs\label{data v1}
&g(L)=-\ln\left[\int_0^{E_{\text{mx}}}\rho(E)\exp\bigg(-\int_L\mu(x,E)\dd x\bigg)\dd E\right].
\es
Here $\rho(E)$ denotes the intensity of the beam at energy $E$ multiplied by the efficiency of the X-ray detector at energy $E$. In practice, CT data are normalized (by performing the so-called offset and gain correction), so we assume $\int_0^{E_{\text{mx}}}\rho(E)\dd E=1$. If $\rho(E)\approx \delta(E-E_0)$ for some $E_0$, then \eqref{data v1} does provide line integral data for the function $\mu(x,E_0)$:
\bs\label{data v11}
&g(L)=\int_L\mu(x,E_0)\dd x.
\es
In most cases, the spectrum $\rho(E)$ is spread out; therefore, the assumption $\rho(E)\approx \delta(E-E_0)$ is not accurate. 

Ignoring the dependence of $\mu$ on $E$ during reconstruction leads to the appearance of beam-hardening artifacts (BHA). There are three types of BHA: cupping (see Figure~\ref{fig:global wide}), dark localized strips between objects with an effective atomic number higher than the background, and streaks (see Figure~\ref{fig:global narrow} for the last two artifacts).

For some objects, it is possible to remove the BHA using CT data preprocessing \cite{Joseph1978}, \cite[p. 795]{Schuller2015} and \cite[Section 12.5]{eps08}. Suppose $\mu(x,E)=\mu(E) f(x)$. Here, $\mu$ is a {\it known} function that depends only on energy, and $f(x)$ is an {\it unknown} function that depends only on position and needs to be reconstructed. This approximation is fairly accurate, for example, for biological soft tissues, because they consist primarily of water. Therefore, when imaging anatomical regions with minimal amounts of bone, we obtain $\mu(E)=\mu_w(E)$. The latter is the well-known attenuation coefficient of water as a function of energy. Consequently, the function $f(x)$ is the water-calibrated density of the tissue at point $x$. 

Assuming the above model, we can rewrite \eqref{data v1} as 
\bs\label{data v1_alt}
&g(L)=G(t),\ G(t)=-\ln\left[\int_0^{E_{\text{mx}}}\rho(E)\exp\big(-\mu(E)t\big)\dd E\right],\\ &t=(Rf)(L)=\int_Lf(x)\dd x,
\es
where $R$ denotes the classical two-dimensional (2D) Radon transform. The function $G$ increases monotonically and its inverse, $G^{-1}$, can be easily obtained. Subsequently, $G^{-1}(g(L))$ provides accurate line integral data for $f$. When $\mu=\mu_w$, this procedure is called soft tissue correction \cite{Joseph1978} or water precorrection \cite[p. 795]{Schuller2015}.

However, most of the time, the assumption $\mu(x,E)=\mu(E) f(x)$ is not accurate. For example, bones and soft tissues/water have very different dependencies of the attenuation coefficient on energy in the interval $[0,E_{\text{mx}}]$ \cite[Fig. 1]{Joseph1978}, \cite{kmpk10, Wang2011}. Similarly, metals (such as pacemakers and hip implants), which may appear inside the human body, exhibit a very different dependence of the attenuation coefficient on energy from water and bone.

As a more accurate approximation, one can suppose that 
\be\label{main model}
\mu(x,E)=\mu_1(E)f_1(x)+\mu_2(E)f_2(x). 
\ee
Here, $\mu_k(E)\ge 0$, $k=1,2$, represent the dependence of the two basis attenuation functions on the energy of the X-rays, and $f_k$, $k=1,2$, are density-like functions. The latter reflect some intrinsic properties of the medium to which X-rays are sensitive and which are independent of energy. Typically, they reflect the mass density (or electron density) and effective atomic number $Z_{\text{eff}}$ of the medium. 

It turns out that model \eqref{main model} is highly accurate in most applications (both medical and nonmedical). The model was originally proposed in \cite{Alvarez1975} and has since become known as the Alvarez–Macovski dual-energy model. After its publication, the model became widely accepted. See \cite{Wang2011, Alvarez2011a, Paziresh2016, Sellerer2019, McCollough2020} and the references therein, which are a small sample of the large body of literature that utilizes the model. At present, \eqref{main model} is the only widely accepted model to account for the dependence of the X-ray attenuation coefficient $\mu$ on energy, which is used by practitioners for image reconstruction in the vast majority of cases. 

In some cases, \eqref{main model} is modified to include more terms $\mu_k(E)f_k(x)$. For example, a patient undergoing CT may be injected with a contrast agent that contains high atomic number (high-$Z$) substances such as iodine ($Z = 53$) and/or gadolinium ($Z = 64$). The attenuation coefficients $\mu(E)$ of these elements exhibit a sharp discontinuity (known as the K-edge discontinuity) at energies relevant to medical X-ray imaging. The presence of the K-edge in the data requires the addition of more terms to the model \eqref{main model} \cite{Roessl2007, Alessio2013, McCollough2020}. In this paper we only consider \eqref{main model} with two terms, but our methods can be easily adapted to models with any number of basis functions (see Remark~\ref{rem:bdries} and Remark~\ref{rem:num basis}).

Reconstructing an accurate artifact-free image from data \eqref{data v1} is difficult. This follows from the fact that the energy dependence of $\mu(x,E)$ is averaged over all energies in the X-ray beam. Reducing BHA requires complicated processing. Customized algorithms for BHA correction have been developed; however, they have drawbacks and are only guaranteed to work under specific circumstances \cite{Dewulf2012, Schuller2015, mfda18}. Frequently, beam hardening is ignored, and it is assumed that the reconstructed image represents the attenuation coefficient at an effective energy, $\mu(x,E_{\text{eff}})$. 

Consequently, the fact that CT data are not line integrals represents the following important question that is not yet fully resolved. What is the nature of the reconstructed image and, more precisely, of image artifacts when conventional Radon transform inversion is applied to the data \eqref{data v1}?   

Two early papers containing a semi-heuristic analysis of nonlinear effects in CT are \cite{Palam86, Palam90}. An important step towards rigorously answering this question is made in \cite{Park2017a}. The authors propose a model based on the following approximations:
\begin{enumerate}
\item The source spectrum is confined to a narrow interval $[E_0-\e,E_0+\e]$.
\item The materials in the object are classified as metal and non-metals. 
\item The attenuation coefficient of the non-metals is approximately constant in this interval, while the attenuation coefficient of the metal is approximately linear in the interval. 
\end{enumerate}
Hence, the attenuation coefficient is approximated by
\be\label{alt mu}
\mu(x,E)=\mu(x,E_0)+\al(E-E_0)\chi_D,\ |E-E_0|\le \e,
\ee
where $D$ is the support of the metal object and $\chi_D$ is the characteristic function of $D$. Under the above assumptions, the following nonlinear model for CT data $g$ is derived in \cite{Park2017a}:
\bs\label{alt model}
g=&-\ln\left[\int_{E_0-\e}^{E_0+\e}\rho(E)\exp\big(-R \mu(\cdot,E)\dd x\big)\dd E\right]\\ 
\approx& R \mu(\cdot,E_0)-\ln\bigg(\frac{\sinh(\al\e R\chi_D)}{\al\e R\chi_D}\bigg).
\es
The authors also perform a microlocal analysis of artifacts caused by nonlinearities in the data in the framework of their model. Subsequent microlocal  analyses of the BHA, namely streaks, under the assumption of the same model are in \cite{Palacios2018a, Wang2021a, chih22, Wang2023}. The main finding in these papers is that the data $g$ is locally a paired Lagrangian distribution, and this is the root cause of the streaks in the reconstructed image. 

%While the results obtained by analyzing the model \eqref{alt model} are valuable, their utility is limited by the fact that the model is inaccurate. As the preceding discussion shows, assumptions (1) and (3) are not satisfied in the majority of CT scans.
%Apart from the theoretical works \cite{Palacios2018a, Wang2021a, Wang2023}, the model does not seem to be used in practice.

While the results obtained by analyzing the model \eqref{alt model} are valuable, their utility is limited by the fact that the model's underlying assumptions may not hold. As the preceding discussion shows, assumptions (1) and (3) are not satisfied in the majority of CT scans.

In this paper we study streak artifacts in CT by using the original, most accurate and widely used by practitioners model \eqref{main model}.  
Assuming the model \eqref{main model}, the so-called after logs CT data, $g(\fp)$, $\fp=(\al,p)$, are
\bs\label{data}
&g(\fp)=-\ln G(\fp),\ \fp\in M:=S^1\times\br,\\
&G(\fp):=\int_0^{E_{\text{mx}}}\rho(E)\exp\big(-[\mu_1(E)\hat f_1(\fp)+\mu_2(E)\hat f_2(\fp)]\big)\dd E.
\es
Here $\hat f_k=R f_k$ is the classical Radon transform (CRT) of $f_k$, $k=1,2$. Clearly, $g(\al,p)=g(\al+\pi,-p)$.

Even though \eqref{data} is quite different from \eqref{alt model}, there are major similarities in the analyses of the two models and in the main conclusions. We also show that the root cause of the streaks is data nonlinearity, which renders $g$ a paired Lagrangian distribution. Hence, we can use some of the intermediate results obtained in \cite{Palacios2018a, Wang2021a, Wang2023}. Nevertheless, our approach differs from these papers. There, the authors (1) approximate the analytic functions $\ln(t)$ and $\sinh(t)$ in \eqref{alt model} by polynomials, and then (2) study the singularities of $g$ by analyzing the {\it approximate} map $\mu\to g$ using the notion of paired Lagrangian distributions. 

Instead, we characterize the singularities of $g$ using the approach developed in \cite{rz1, rz5, airapramm01}. This allows us to study the singularities of the function, $\fr$, reconstructed from $g$ in a much easier fashion. In particular, no approximation of the model in \eqref{data} is required. 
Besides using a more accurate model and a different method of analysis, we also explicitly compute the leading singularity of $\fr$ near the beam hardening streaks. This allows us to characterize streaks more precisely than by computing the order of the distribution. A somewhat similar approach was used in \cite{Wang2023} for a related problem.

Interestingly, the model \eqref{data} has already been used for rigorous analysis of tomographic image reconstruction. See \cite{Bal2020}, where the authors prove the uniqueness of reconstruction from dual energy data.

Let $\s$ be a piecewise smooth curve without self-intersections with finitely many disconnected components. We assume that $\ssp(f_k)\subset\s$, $k=1,2$.
The remainder of this paper is organized as follows. In section~\ref{sec:idea} we discuss the problem setting. In section~\ref{sec:main res} we formulate our main result, Theorem~\ref{thm:main}. We also describe qualitatively the singularities of $\fr$ that appear because of data nonlinearity. The proof of Theorem~\ref{thm:main} is spread over several sections. 
In section~\ref{sec:sing g} we describe the singularities of $g$ corresponding to smooth sections of $\s$. In section~\ref{sec:corner} we describe the singularities of $g$ corresponding to a corner of $\s$. In section~\ref{sec:proof} we finish the proof of Theorem~\ref{thm:main}. In sections~\ref{sec:pp} and \ref{sec:pm} we obtain the leading singular behavior of $\fr$ in a neighborhood of streaks caused by data nonlinearity. Finally, in section~\ref{sec:numerics} we conduct a numerical experiment to validate our description of beam hardening streaks in $\fr$.

\section{Problem setting}\label{sec:idea}

\subsection{Assumptions} \label{sec:assump}
We begin by describing our assumptions about the $f_k$.
Suppose there are $J\in\BN$ bounded domains, $D_j\subset\br^2$, $j=1,2,\dots,J$, in the plane. Their boundaries, $S_j:=\pa D_j$, are piecewise smooth and disjoint, $S_j\cap S_k=\varnothing$, $j\not=k$. No assumption about the convexity of the $D_j$’s is made. Denote $\s:=\cup_{j=1}^J \s_j$. If $\s$ is nonsmooth at $x$, then $x$ is a corner point of $\s$ and the two one-sided tangents to $\s$ at $x$ do not coincide. Given $x\in\s$ where $\s$ is smooth, $\varkappa(x)$ denotes the curvature of $\s$ at $x$. 

\begin{definition}\label{def:gen tang}
A line $L$ is tangent to $\s$ at $x$ in the {\it generalized} sense if $L$ is tangent to $\s$ in the conventional sense if $\s$ is smooth at $x$ or contains $x$ if $x$ is a corner of $\s$.
\end{definition}

\begin{definition}\label{def:excp line}
A line $L$ is called {\it exceptional} if it is tangent to $\s$ at a point $x$ where $\varkappa(x)=0$ or if it is a one-sided tangent at a corner point of $\s$. 
\end{definition}

We suppose that (1) no exceptional line is tangent to $\s$ in the generalized sense at another point, and (2) no line is tangent to $\s$ in the generalized sense at more than two points.

Finally, we suppose there are functions $f_{kj}\in C^\infty(\br^2)$ such that 
\be\label{fk model}
f_k=\sum_{j=1}^J f_{k,j}\chi_j,\ k=1,2,
\ee
where $\chi_j$ denotes the characteristic function of $D_j$. Therefore, $WF(f_k)\subset N^*\s$.

\begin{remark}\label{rem:bdries} The model~\eqref{fk model} reflects the common assumption that the object being scanned can be segmented into distinct regions $D_j$, such that the material composition of the object varies smoothly within each region and may change abruptly from one region to another. Hence it is reasonable to assume that, regardless of the number of basis functions, all the $f_k$’s are smooth within each region and may be discontinuous across region boundaries. Therefore, the assumption that $\ssp f_k\subset \s$ for all $k$ is practically reasonable.
\end{remark}

Throughout the paper, we use the following convention. Let $I$ be a space of distributions, for example, conormal or paired Lagrangian distributions. Recall that $\fp=(\al,p)$ and $M=S^1\times\br$ (see \eqref{data}). We say that $g\in I$ {\it near $\fp_0$ (or in an open set $M_0$)} when $\chi g\in I$ for any $\chi\in C^\infty(M)$ supported in a sufficiently small neighborhood of $\fp_0$ (or supported  in $M_0$). Likewise, given a set $\Xi\subset T^*M$, we say that $WF(g)\in \Xi$ {\it near $\fp_0$ (or, in an open set $M_0$)} if $WF(\chi g)\in \Xi$ for any $\chi$ as above.

\subsection{Spaces of symbols and distributions}

Let $U\subset\br^n$ be an open set. We say that $a\in S^r(U\times\br^N)$ is the standard symbol of order $r\in\br$ if for any compact $K\subset U$,
\be
\big|\pa_x^\bt\pa_\xi^\ga a(x,\xi)\big|\le C_{\bt,\ga,K} (1+|\xi|)^{r-|\ga|},\ 
x\in K,\xi\in\br^N,\bt\in \BN_0^n,\ga\in \BN_0^N.
\ee
Here 
%$\langle\xi\rangle=(1+|\xi|^2)^{1/2}$ and 
$\BN_0=0\cup\BN$. Let $\Lambda\subset T^*U$ be a closed conic Lagrangian submanifold. Suppose $\Lambda$ is parameterized by a homogeneous phase function $\phi:U\times\br^N\to\br$:
\be
\Lambda=\{(x,\dd_x\phi(x,\xi))\in T^*U: \dd_\xi\phi(x,\xi)=0\}.
\ee
Then the oscillatory integral
\be\label{lagr distr}
u(x)=\int_{\br^N}a(x,\xi) e^{i\phi(x,\xi)}\dd\xi,\ x\in U,
\ee
where $a\in S^r(U\times\br^N)$ defines a Lagrangian distribution $u\in I^q(\Lambda)$ of order $q=r-(n/4)+(N/2)$. If $\s\subset U$ is a smooth submanifold and $\Lambda=N^*\s$ is the conormal bundle of $\s$, then $u$ is a conormal distribution: $u\in I^q(U;\s)$ \cite[Section 18.2]{hor3} and \cite[Section 25.1]{hor4}.

Suppose $U=U_1\times U_2$, where $U_j\subset\br^{n_j}$,  $x=(x_1,x_2)$, $x_j\in\br^{n_j}$, $j=1,2$, and $n=n_1+n_2$. Then a Lagrangian distribution $u\in I^q(\Lambda)$ defines a Fourier Integral Operator (FIO) $\CU:\coi(U_2)\to \CD^\prime(U_1)$, $\CU\in I^q(U_1\times U_2,\Lambda)$, of order $q$ by the formula
\be
(\CU f)(x_1)=\int_{\br^N}\int_{U_2}a(x,\xi) e^{i\phi(x,\xi)}f(x_2)\dd x_2\dd\xi,\ x_1\in U_1.
\ee
See \cite[Section 25.2]{hor4}.

Suppose that two Lagrangian submanifolds $\Lambda_0,\Lambda_1\subset T^*U$ intersect trans\-versally at a codimension $k$ submanifold $\Omega$. The space of paired Lagrangian distributions of order $p,l$, denoted $I^{p,l}(\Lambda_0,\Lambda_1)$, is the set of distributions $u$ that satisfy $u\in I^{p+l}(\Lambda_0\setminus\Omega)$ and $u\in I^p(\Lambda_1\setminus\Omega)$, see \cite{meluhl, guuhl} for more details.

\subsection{Preliminaries} 
Pick any $x_0\in\s$. Let $U$ be a small neighborhood of $x_0$. As is well known, $f_k\in I^{-1}(U;\s)$ if $\s\cap U$ is smooth. If $x_0$ is a corner of $\s$, we assume $\s\cap U=\s_1\cup\s_2$, where $\s_{1,2}$ are smooth curve segments that share a common endpoint, $x_0$ (see Figure~\ref{fig:corner}, left panel). In this case 
\be
f_k\in I^{-1,-1/2}(N^* x_0;\s_1)+I^{-1,-1/2}(N^* x_0;\s_2)
\ee
near $x_0$ \cite[Lemma 5.1]{Palacios2018a}. 

The CRT $R:\mathcal E(\br^2)\to\CD(M)$ is an elliptic FIO given by
\be\label{rt fio}
\hat f(\fp)=\int_{\br^2}f(x)\de(\vec\al\cdot x-p)\dd x
=\frac1{2\pi}\int_{\br}\int_{\br^2}f(x)e^{i\la(\vec\al\cdot x-p)}\dd x\dd\la.
\ee
This implies that the CRT is an FIO of order $q=0-(n/4)+(k/2)=-1/2$, where $n=4$ and $k=1$. The ambient space here is $M\times\br^2$, and thus, its dimension is $n=4$. The canonical relation of the CRT is
\bs\label{crt canon rel}
C=\big\{&\big((\al,p),\la(\vec\al^\perp\cdot x,-1));(x,-\la\vec\al)\big):\\
&(\al,p)\in M,x\in\br^2,\vec\al\cdot x-p=0,\la\in\br\setminus 0 \big\}\subset T^*M\setminus 0\times T^*\br^2\setminus0.
\es
Clearly, $C$ is a homogeneous canonical relation. A CRT inversion formula is as follows:
\be\label{rad inv orig}
f(x)=(R^{-1}\hat f)(x)=-\frac1{2\pi^2}\int_0^{\pi}\int_{\br} \frac{\pa_p g(\al,p)}{p-\vec\al\cdot x}\dd p\dd\al=-\frac1{2\pi}(R^*\CH \pa_p f)(x),
\ee
where $\CH$ denotes the Hilbert transform. This is also an elliptic FIO, and its homogeneous canonical relation, $C^*$, is obtained from $C$ in \eqref{crt canon rel} by swapping the coordinates:
\bs\label{crtstar canon rel}
C^*=\big\{&\big((x,\xi);(\fp,\eta)\big):\big((\fp,\eta);(x,\xi)\big)\in C \big\}\subset T^*\br^2\setminus0\times T^*M\setminus 0.
\es

It is clear from \eqref{crt canon rel} and \eqref{data} that 
\bs\label{all sings}
&\ssp(\hat f_k),\ssp(g)\subset \Gamma,\\
&\Gamma=\{(\al,p)\in M:\ (x,\vec\al)\in N^* \s,\ p=\vec\al\cdot x\}.
\es
Put simply, $\fp\in\Gamma$ if $L_{\fp}:=\{x\in\br^2:\vec\al\cdot x=p\}$ is tangent to $\s$ in the generalized sense. Thus, $\Gamma$ is the generalized Legendre transform of $\s$ \cite{rz5}. We denote this relation by $\Gamma=\CL(\s)$. The generalized Legendre transform is very closely related to the Legendre-Fenchel transform \cite[Chapter VI]{Ellis85}.

\section{Main results}\label{sec:main res}

We now state our main results.

\begin{theorem}\label{thm:main}
Let $f_{1,2}$ and $\s$ satisfy the assumptions in section~\ref{sec:assump}.
Select some $\fp_0\in \Gamma$ and let $M_0$ be a sufficiently small neighborhood of $\fp_0$. One has
\begin{enumerate}
\item\label{clssp} $\ssp (g)\subset \Gamma$.
\item\label{clnonsm} If $\Gamma$ is nonsmooth at $\fp_0$, then $L_{\fp_0}$ is tangent to $\s$ at an isolated point where the curvature of $\s$ equals zero. 
\item\label{clselfint} If $\Gamma$ self-intersects at $\fp_0$, i.e. $\Gamma\cap M_0=\Gamma_1\cup\Gamma_2$ for some smooth $\Gamma_{1,2}$, then $L_{\fp_0}$ 
\begin{enumerate}
\item is tangent to $\s$ at two different points, or 
\item contains a line segment of $\s$, or 
\item is a one-sided tangent to $\s$ at a corner point. 
\end{enumerate}
In cases (a) and (b), $\Gamma_1$ and $\Gamma_2$ are transversal at $\fp_0$, and $\Gamma_1$ and $\Gamma_2$ are tangent at $\fp_0$ in case (c).
\item\label{clorder} If $\Gamma\cap M_0$ is a smooth curve without self-intersections, then $\hat f_k\in I^{-r}(M_0;\Gamma)$, $k=1,2$, implies $g\in I^{-r}(M_0;\Gamma)$. Here $r=3/2$ if $L_{\fp_0}$ passes through a corner of $\s$ and $r=2$ if $L_{\fp_0}$ is tangent to $\s$.
\item\label{clinters} Suppose $\Gamma\cap M_0=\Gamma_1\cup\Gamma_2$, where $\Gamma_{1,2}$ are smooth and intersect transversely at $\fp_0$. Because of nonlinear interaction, $g$ is a paired Lagrangian distribution near $\fp_0$:
\be
g\in I^{-r_1,-r_2+(1/2)}(N^*\fp_0,N^*\Gamma_1)+I^{-r_2,-r_1+(1/2)}(N^*\fp_0,N^*\Gamma_2).
\ee
Here $r_j=3/2$ in the case of the usual tangency at $x_j$ and $r_j=2$ if $x_j$ is a corner point, $j=1,2$.
\item\label{clrest} Let $M^\prime$ be the finite collection of all $\fp\in M$ such that either $\Gamma$ is nonsmooth at $\fp$ or self-intersects at $\fp$. Then
\be\label{WFg full}
WF(g)\subset N^*(\Gamma\setminus M^\prime)\cup \big(\cup_{\fp\in M^\prime}N^*\fp\big).
\ee
\end{enumerate}
\end{theorem}

Let us briefly discuss the above six claims. The first three are about the geometric properties of the singular support of $g$, $\Gamma$, and how they relate to the geometric properties of $\s$. The last three claims provide a microlocal description of the singularities of $g$.

Claim~\eqref{clssp} states that the singular support of $g$ does not extend beyond $\Gamma$, the singular support of $\hat f_k$'s. Claims~\eqref{clnonsm} and \eqref{clselfint} are well-known properties of the Legendre transform. We state and prove them to make the proofs of the other claims self-contained. Claim~\eqref{clorder} asserts that if $\Gamma\cap M_0$ is smooth, then $g$ and $\hat f_k$ are locally conormal distributions of the same order. 

\begin{figure}[h]
{\centerline{
{\epsfig{file={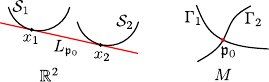}, width=7cm}}
}}
\caption{The case when local segments of $\s$ near $x_1$ and $x_2$ are on one side of the double tangent $L_{\fp_0}$.}
\label{fig:standard}
\end{figure}

Claim~\eqref{clinters} states that if $\Gamma\cap M_0=\Gamma_1\cap\Gamma_2$ and the intersection is transversal, then $g$ is a paired Lagrangian distribution with the given orders. Such a shape of $\Gamma$ occurs near a point $\fp_0$ if $L_{\fp_0}$ is a double tangent to $\s$. See, e.g., Figure~\ref{fig:standard}.

Claims~\eqref{clorder} and \eqref{clinters} do not address the following two cases: (1) $\Gamma$ is nonsmooth at $\fp_0$; and (2) $\Gamma\cap M_0=\Gamma_1\cap\Gamma_2$ and $\Gamma_{1,2}$ are tangent to each other at $\fp_0$. Claim~\eqref{clrest} states that, at worst, this singular behavior of $\Gamma$ at $\fp_0\in M^\prime$ adds $N^*\fp_0$ to $WF(g)$. We do not study these two cases in more detail because they are not related to data nonlinearity.

%\subsection{Singularities of $\fr$}\label{ssec:sing_f}

Since $R^{-1}$ is an FIO with the canonical relation given in \eqref{crtstar canon rel}, \eqref{WFg full} gives
\be\label{WFf full}
WF(f^\text{rec})\subset N^*\s\cup \big(\cup_{\fp\in M^\prime}N^*L_{\fp}\big).
\ee
In other words, whenever $N^*\fp$ appears in $WF(g)$, it may cause a streak along $L_\fp$ in $\fr$. The reason $N^*\fp$ appears in $WF(g)$ (or does not appear) can be different. To better understand the implications of \eqref{WFf full}, we identify the following three groups of cases. 

\begin{enumerate}
\item $\Gamma$ is locally smooth (claim~\eqref{clorder} of Theorem~\ref{thm:main}). In this case, no singularities are added to $\fr$.
\item Locally, $\Gamma=\Gamma_1\cap\Gamma_2$, where $\Gamma_j$ are smooth and intersect transversely at $\fp_0$ (claim~\eqref{clinters} of Theorem~\ref{thm:main}). In this case, the component $N^*\fp_0\subset WF(g)$ appears owing to {\it nonlinear interaction}. This can cause a streak along $L_{\fp_0}$ in $\fr$.
\item Exceptional cases: $\Gamma$ is not smooth or $\Gamma=\Gamma_1\cap\Gamma_2$, where $\Gamma_j$ are smooth and tangent to each other at some $\fp\in M$ (all the remaining cases in claim~\eqref{clrest} of Theorem~\ref{thm:main}). The component $N^*\fp\subset WF(g)$ also appears in this case, leading to a possible streak along $L_{\fp}$ in $\fr$. However, the root cause of the streak is not the nonlinearity of the data measurement, but the fact that $g$ is not in the range of $R$.
\end{enumerate}

See also \cite{Palam86, Palam90}, where cases (2) and (3), which describe artifacts, are identified using heuristic arguments.

%The leading singularity of $f$ in a neighborhood of artifacts caused by nonlinearity is obtained in the next two sections.

\section{Singularities of data -- locally smooth $\s$}\label{sec:sing g}
Our first goal is to describe the behavior of $g(\fp)$ near its singular support. Select some $\fp_0=(\al_0,p_0)\in\Gamma$ and consider the line $L_{\fp_0}$. Suppose $L_{\fp_0}$ is tangent to $\s$ at some point $x_0$, $\varkappa(x_0)\not=0$. Let $p=P(\al)$ be a parameterization of $\Gamma$ in a neighborhood of $\fp_0$. The assumption $\varkappa(x_0)\not=0$ implies that both $\Gamma$ and $P(\al)$ are locally smooth. It is easy to establish that
\bs\label{sing rt}
&\hat f_{k}(\fp)=(p-P(\al))_{\iota}^{1/2}R_{k}(\fp)+T_{k}(\fp),\ \fp\in M_0,\\
&R_{k}(\fp_0)=2(2/\varkappa(x_0))^{1/2},
\es
where $\iota=+$ or $-$ (depending on the location of $\s$ relative to $L_{\fp_0}$), $t_\iota=\max(0,\iota t)$, $M_0=I_\al\times I_p$ is a sufficiently small, rectangular neighborhood of $\fp_0$, $R_{k},T_{k}\in C^\infty(M_0)$, 
%$\varkappa_j(x_0)$ is the curvature of $\s_j$ at $x_0$, 
and $P\in C^\infty(I_\al)$ \cite{rz1, rz5, airapramm01}.

\subsection{Tangency at a single point}\label{ssec:single}
Suppose $L_{\fp_0}$ is tangent to $\s$ {\it only at one point}, $x_0$, and $L_{\fp_0}$ is not tangent to $\s$ in the generalized sense at another point. We substitute \eqref{sing rt} into \eqref{data}. Expand $e^t$ at $t=0$ in the Taylor series and group the terms to get
\bs\label{exp exp}
G(\fp)=&(p-P(\al))_{\iota}^{1/2}A_{1/2}(\fp)+(p-P(\al))_{\iota}A_1(\fp)
+G_{0}(\fp),\ \fp\in M_0,
\es
where $A_{1/2},A_1,G_{0}\in C^\infty(M_0)$. These functions are given by
\bs\label{exp exp aux}
A_{1/2}(\fp)=&-\sum_{n\ge 1,n\text{ odd}}\frac{G_n(\fp)}{n!}(p-P(\al))^{(n-1)/2},\\
A_1(\fp)=&\sum_{n\ge 2,n\text{ even}}\frac{G_n(\fp)}{n!}(p-P(\al))^{(n-2)/2},\\
G_n(\fp)=&\int _0^{E_{\text{mx}}}W(E,\fp)\big[\mu_1(E)R_{1}(\fp)+\mu_2(E)R_{2}(\fp)\big]^n\dd E,\\
W(E,\fp)=&\rho(E) \exp\big(-[\mu_1(E)T_{1}(\fp)+\mu_2(E)T_{2}(\fp)]\big).
\es
The two series converge absolutely because $|G_n(\fp)|\le c^n$, $\fp\in M_0$, for some $c>0$.

Substituting \eqref{exp exp} into the first equation in \eqref{data} yields
\bs\label{g st}
g=&-\ln\big(G_0+(p-P(\al))_\iota^{1/2}A_{1/2}+(p-P(\al))_\iota A_1\big)\\
=&-\ln G_0+\sum_{n\ge 1}\frac{(-1)^n}{n}\frac{\big[(p-P(\al))_\iota^{1/2}A_{1/2}+(p-P(\al))_\iota A_1\big]^n}{G_0^n},
\es
where we dropped $\fp$ from the arguments of all functions. We assume that there is no photon starvation, that is,
\be
\inf_{\fp\in M_0}\Vert W(\cdot,\fp)\Vert_{L^\infty([0,E_{\text{mx}}])}>0.
\ee
In this case, $\inf_{\fp\in M_0}G_0(\fp)>0$. Clearly, the series in \eqref{g st} converges absolutely provided $\iota(p-P(\al))>0$ is sufficiently small. If $\iota(p-P(\al))\le0$, then the series vanishes. Therefore,
\be\label{g sing tan}
g(\fp)=-\ln G_0(\fp)+(p-P(\al))_\iota^{1/2}B_{1/2}(\fp)+(p-P(\al))_\iota B_1(\fp),\ \fp\in M_0,
\ee
for some $B_r\in C^\infty(M_0)$, $r\in\{1/2,1\}$. In particular, by \eqref{exp exp aux}, 
\bs
B_{1/2}=&\frac{G_1}{G_0}=\frac{\int _0^{E_{\text{mx}}}W(E) \big[\mu_1(E)R_{1}+\mu_2(E)R_{2}\big]\dd E}{\int _0^{E_{\text{mx}}}W(E)\dd E},\ p=P(\al),\al\in I_\al.
\es
As is easily seen from \eqref{sing rt} and \eqref{g sing tan}, 
\be\label{conorm tang}
\hat f_k,g\in I^{-3/2}(M_0;\Gamma)
\ee
are conormal distributions. The claim regarding $\hat f_k$ is well known. Since 
\be
\CF(t_\pm^{1/2})=O(|\la|^{-3/2}),\ \la\to\infty, 
\ee
(see \eqref{p+} and \eqref{pminus}), the order of the conormal distributions $\hat f_k$ and $g$ is $q=(-3/2)+(n/4)-(k/2)=-3/2$. Here, $n=2$ is the dimension of the ambient space $\br^2$, and $k=1$ is the dimension of the dual variable $\la$.

%Assuming that there is only one material and the X-ray beam is monoenergetic, i.e. $\rho(E)=\de(E-E_0)$ for some $E_0\in (0,E_{mx})$, we recover the well-known formulas. \red{What happens to $(p-P_j)_+$?}

\subsection{Tangency at two points}\label{ssec:double}
Next, suppose $L_{\fp_0}$ is tangent to $\s$ at two points $x_1\not=x_2$ (see Figure~\ref{fig:standard}). Let $\s_j$ be a local segment of $\s$ in a neighborhood of $x_j$, $j=1,2$. We assume that $\s_{1,2}$ are smooth and $\varkappa(x_j)\not=0$, $j=1,2$. Further, $\Gamma\cap M_0=\Gamma_1\cup\Gamma_2$, where $\Gamma_j=\CL(\s_j)$. Additionally, $\Gamma_1$ and $\Gamma_2$ intersect transversely at $\fp_0$ \cite{Palacios2018a}. Let $p=P_j(\al)$ be a local parameterization of $\Gamma_j$. Then $p_0=P_1(\al_0)=P_2(\al_0)$. The analog of \eqref{sing rt} becomes
\bs\label{sing rt v2}
&\hat f_{k}(\fp)=(p-P_j(\al))_{\iota_j}^{1/2}R_{k,j}(\fp)+T_{k,j}(\fp),\ \fp\in M_0,\\
&R_{k,j}(\fp_0)=2(2/\varkappa(x_j))^{1/2},\ j=1,2.
\es
Substituting \eqref{sing rt v2} into \eqref{data} gives
\be
G(\fp)=\sum_{m,n\ge0}\frac{(-1)^{m+n}}{m!n!} G^{(m,n)}(\fp),\ \fp\in M_0,
\ee
where the series converges absolutely and
\bs\label{aux}
G^{(m,n)}:=&\int _0^{E_{\text{mx}}}W(E)\big[(p-P_1)_{\iota_1}^{1/2}\big(\mu_1(E)R_{1,1}+\mu_2(E)R_{2,1}\big)\big]^m\\
&\qquad \times \big[(p-P_2)_{\iota_2}^{1/2}\big(\mu_1(E)R_{1,2}+\mu_2(E)R_{2,2}\big)\big]^n\dd E,\\
W(E)=&\rho(E) \exp\big[-\mu_1(E)(T_{1,1}+T_{1,2})-\mu_2(E)(T_{2,1}+T_{2,2})\big].
\es
For simplicity, we removed the arguments $\al$ and $p$ from all the functions in the above equation. Expanding further and combining terms gives
\bs\label{AtoG}
G(\fp)=&A_{0,0}(\fp)+\sum_{s_1\in\{1/2,1\}}G_{s_1,0}(\fp)+\sum_{s_2\in\{1/2,1\}}G_{0,s_2}(\fp)\\
&+\sum_{s_1,s_2\in\{1/2,1\}}G_{s_1,s_2}(\fp),\ \fp\in M_0,\\
G_{s_1,s_2}(\fp):=&(p-P_1(\al))_{\iota_1}^{s_1} (p-P_2(\al))_{\iota_2}^{s_2} A_{s_1,s_2}(\fp),\\
A_{0,0}(\fp)=&\int _0^{E_{\text{mx}}}W(E,\fp)\dd E,
\es
where $A_{s_1,s_2}\in C^\infty(M_0)$.

As in the preceding subsection, if there is no photon starvation, then 
\be
\inf_{\fp\in M_0}A_{0,0}(\fp)>0. 
\ee
Expand $\ln t$ in a series at $t=1$ to get 
\bs\label{g exp v2}
g(\fp)=B_{0,0}(\fp)+&\sum_{s_1\in\{1/2,1\}}(p-P_1(\al))_{\iota_1}^{s_1} B_{s_1,0}(\fp)\\
+&\sum_{s_2\in\{1/2,1\}}(p-P_2(\al))_{\iota_2}^{s_2} B_{0,s_2}(\fp)\\
+&\sum_{s_1,s_2\in\{1/2,1\}}(p-P_1(\al))_{\iota_1}^{s_1} (p-P_2(\al))_{\iota_2}^{s_2} B_{s_1,s_2}(\fp),\\
B_{0,0}(\fp)=&-\ln\big(A_{0,0}(\fp)\big),\ \fp\in M_0,
\es
where $B_{s_1,s_2}\in C^\infty(M_0)$, $s_1,s_2\in\{0,1/2,1\}$. In particular, setting $\fp=\fp_0$ gives
\bs\label{B main}
&B_{\frac 12,\frac12}=A_{\frac12,\frac12}/A_{0,0}\\
&=\frac{\int _0^{E_{\text{mx}}}W(E) \big[\mu_1(E)R_{1,1}+\mu_2(E)R_{2,1}\big]\big[\mu_1(E)R_{1,2}+\mu_2(E)R_{2,2}\big]\dd E}{\int _0^{E_{\text{mx}}}W(E)\dd E}.
\es

By \eqref{g exp v2},
\be\label{tang pair lagr}
g\in I^{-3/2,-1}(N^*\fp_0,N^*\Gamma_1)+I^{-3/2,-1}(N^*\fp_0,N^*\Gamma_2)
\ee
near $\fp_0$. 

To explain the values $p=-3/2$ and $l=-1$ for the superscripts of $I^{p,l}$ in \eqref{tang pair lagr}, note that $g\in I^{-3/2}(N^*\Gamma_j\setminus N^*\fp_0)$, $j=1,2$, by \eqref{conorm tang}. Therefore $p=-3/2$. The value of $l$ is computed from the requirement $g\in I^{-3/2}(N^*\fp_0\setminus N^*\Gamma_j)$. Hence, we solve $p+l=-3-(n/4)+(k/2)$. The number -3 in this equation appears because each factor $(p-P_j(\al))_{\iota_j}^{1/2}$ contributes a decay of order $-3/2$ in the Fourier domain, and we take $k=2$ because this is the dimension of the Fourier integral representing the product $(p-P_1(\al))_{\iota_1}^{1/2}(p-P_2(\al))_{\iota_2}^{1/2}$.

\begin{remark}\label{rem:num basis} Equation \eqref{g exp v2} also holds for models of the kind \eqref{main model} with more basis functions $\mu_k(E)f_k(x)$, $k=1,2,\dots,K$, as long as $\ssp f_k\in\s$ for all $k$. The only difference is that the formulas for the coefficients $A_{s_1,s_2}$ involve more terms (cf. \eqref{aux} and  \eqref{AtoG}). 
\end{remark}

\section{Corner cases.}\label{sec:corner}

Suppose $\s$ has a corner at $x_0$, i.e. one-sided tangents to $\s$ at $x_0$, $L_{\fp_1}$ and $L_{\fp_2}$, do not coincide (see Figure~\ref{fig:corner}). Pick $\fp_0\in M$, $\fp_0\not=\fp_{1,2}$, such that $x_0\in L_{\fp_0}$. As before, let $M_0$ be a sufficiently small neighborhood of $\fp_0$. It is easily seen that, if $L_{\fp_0}$ is not tangent to $\s$ elsewhere in the generalized sense, then
\bs\label{sing rt corner}
&\hat f_k(\fp)=(p-P(\al))_+R_k^{+}(\fp)+(p-P(\al))_-R_k^{-}(\fp)+T_k(\fp),\\ 
&P(\al)=\vec\al\cdot x_0,\ \fp\in M_0,
\es
where $R_k^{\pm},T_k\in C^\infty(M_0)$, $k=1,2$. Similarly to section~\ref{ssec:single} we find
\be\label{sing g corner}
g(\fp)=-\ln G_0(\fp)+(p-P(\al))_+B_+(\fp)+(p-P(\al))_-B_-(\fp),\ \fp\in M_0,
\ee
for some $G_0,B_\pm\in C^\infty(M_0)$. Clearly, 
\be\label{conorm corner}
\hat f_k,g\in I^{-2}(M_0;\Gamma_{x_0}).
\ee
Here and throughout the paper we denote
\be 
\Gamma_x=\{(\al,\vec\al\cdot x)\in M:\ \vec\al\in S^1\}.
\ee
Similarly to \eqref{conorm tang}, the order -2 follows from $\CF(t_\pm)=O(|\la|^{-2}),\ \la\to\infty$.

Now we consider the general case when $L_{\fp_0}$ is tangent to $\s$ in the generalized sense at two points, $x_1$ and $x_2$. See, for example, Figure~\ref{fig:two_corner}, where the case of two corner points is shown. Then, $\Gamma\cap M_0=\Gamma_1\cup\Gamma_2$, $\Gamma_1$ and $\Gamma_2$ are smooth and intersect transversely at $\fp_0$. The complete expression for $g$ is too long; therefore, we describe only the relevant paired Lagrangian part:
\bs\label{g exp v3}
g(\fp)-&\sum_{\iota_1,\iota_2\in\{+,-\}}\sum_{s_1,s_2\in\{1/2,1\}}(p-P_1(\al))_{\iota_1}^{s_1} (p-P_2(\al))_{\iota_2}^{s_2} B^{s_1,s_2}_{\iota_1,\iota_2}(\fp)\\
\in & I^{-r_1}(M_0;\Gamma_1)+I^{-r_2}(M_0;\Gamma_2)
\es
with some $B^{s_1,s_2}_{\iota_1,\iota_2}\in C^\infty(M_0)$. The values of $r_j$ and $s_j$ depend on the singularity type of $f$ near $x_j$, $j=1,2$. If $\s$ is smooth at $x_j$, then $r_j=3/2$ and $s_j\in\{1/2,1\}$, cf. \eqref{g sing tan} and \eqref{g exp v2}. If $\s$ has a corner at $x_j$, then $r_j=2$ and $s_j=1$, cf. \eqref{sing g corner}. Furthermore,
\be\label{gen dbl case}
g\in I^{-r_1,-r_2+(1/2)}(N^*\fp_0,N^*\Gamma_1)+I^{-r_2,-r_1+(1/2)}(N^*\fp_0,N^*\Gamma_2)
\ee
in $M_0$. The orders $p$ and $l$ are computed similarly to \eqref{tang pair lagr} using \eqref{conorm tang} and \eqref{conorm corner}.

%\begin{proof} 

\section{End of proof of Theorem~\ref{thm:main}}\label{sec:proof}

Claim~\eqref{clssp} is proven in \eqref{all sings}. To prove claims~\eqref{clnonsm}--\eqref{clrest}, we review several cases. 

%Using the results in \cite{Palacios2018a}, we obtain the following theorem.

\noindent
{\bf 1.} $\s$ is smooth at $x_0$, $\varkappa(x_0)\not=0$, and $L_{\fp_0}$ is tangent to $\s$ only at $x_0$. As is well known, in this case $\Gamma\cap M_0$ is smooth and $\hat f_k\in I^{-3/2}(N^*\Gamma)$ is a conormal distribution near $\fp_0$, $k=1,2$. From \eqref{g sing tan}, $g$ is also a conormal distribution, and $g\in I^{-3/2}(N^*\Gamma)$ near $\fp_0$. 
%In this case no streak artifacts arise in $\fr$ along $L_{\fp_0}$. 

\noindent
{\bf 2.} $L_{\fp_0}$ is tangent to $\s$ at $x_1$ and $x_2$ (and only at these two points), $\s$ is smooth in a neighborhood of $x_1$, $x_2$, and $\varkappa(x_k)\not=0$, $k=1,2$, see Figure~\ref{fig:standard}.  
Let $\s_k$ be a local segment of $\s$ near $x_k$ and $\Gamma_k:=\CL(\s_k)\cap M_0$, $k=1,2$. Then $\Gamma_1\cap\Gamma_2=\fp_0$ and the intersection is transversal \cite[p. 4920]{Palacios2018a}. 

In the linear Radon transform framework, the singularities at $\Gamma_1$ and $\Gamma_2$ do not interact. However, the nonlinearity of the measurement leads to an interaction at $\fp_0$ (see \eqref{g exp v2}). As a consequence, 
\be\label{WFg tt}
WF(g)\subset N^*\Gamma_1\cup N^*\Gamma_2\cup N^*\fp_0 
\ee
near $\fp_0$. More precisely, \eqref{g exp v2} implies (see also \cite[p. 4920]{Palacios2018a}):
\bs\label{tang paired lagr v2}
g\in & I^{-3/2,-1}\big(N^* \fp_0,N^*\Gamma_1\big)
+I^{-3/2,-1}\big(N^*\fp_0,N^*\Gamma_2\big).
\es
The component $N^*\fp_0$ in $WF(g)$ arises owing to the nonlinear interaction of singularities at $\fp_0$.  

\begin{figure}[h]
{\centerline{
{\epsfig{file={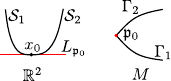}, width=4.5cm}}
}}
\caption{The case $\varkappa(x_0)=0$.}
\label{fig:zero_curv}
\end{figure}

\noindent
{\bf 3.}  $\s$ is smooth near $x_0$, $\varkappa(x_0)=0$ and $L_{\fp_0}$ is tangent to $\s$ at $x_0$, see Figure~\ref{fig:zero_curv}. According to Definition~\ref{def:excp line}, $L_{\fp_0}$ is an exceptional line; therefore, it is not tangent to $\s$ elsewhere in the generalized sense. As is easily seen, $\Gamma$ is non-smooth at $\fp_0$, i.e. $\Gamma\cap M_0 =\Gamma_1\cup\Gamma_2$ for some $\Gamma_1$ and $\Gamma_2$ that are smooth and share an endpoint, $\fp_0$ (see Figure~\ref{fig:zero_curv}, right panel). Then
\be\label{WFg except 1}
WF(\hat f_k),WF(g)\subset N^*\Gamma_1\cup N^*\Gamma_2\cup N^*\fp_0. 
\ee
This implies that the component $N^*\fp_0$ of $WF(g)$ arises because of the shape of $\s$ near $x_0$ and not because of a nonlinear effect. 
%Consequently, a streak artifact in $\fr$ along $L_{\fp_0}$ may appear because $g$ is not in the range of the Radon transform, not because of the nonlinearity of the measurement.

\begin{figure}[h]
{\centerline{
{\epsfig{file={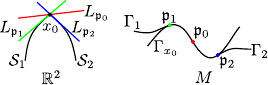}, width=7cm}}
}}
\caption{One corner case.}
\label{fig:corner}
\end{figure}

\noindent
{\bf 4.}  $\s$ has a corner at $x_0$, see Figure~\ref{fig:corner}. Thus, $\s=\s_1\cup\s_2$ in a neighborhood of $x_0$, where $\s_{1,2}$ are smooth curves that share an endpoint, $x_0$. Let $L_{\fp_j}$ be  tangent to $\s_j$ at $x_0$, $j=1,2$. Pick any $\fp_0\in M$ such that $x_0\in L_{\fp_0}$, $\fp_0\not=\fp_{1,2}$, and $L_{\fp_0}$ is not tangent to $\s$ in the generalized sense elsewhere. Because $L_{\fp_1}$ and $L_{\fp_2}$ are one-sided tangents to $\s$ at $x_0$, they are exceptional lines and are not tangent to $\s$ in the generalized sense elsewhere either. 

Let $M_0$ be a small neighborhood of $\{\fp_0,\fp_1,\fp_2\}$ and denote $\Gamma_j=\CL(\s_j)\cap M_0$, $j=1,2$.
%Denote 
%\be 
%\Gamma_{x}=\{(\al,\vec\al\cdot x):\ \vec\al\in S^1\},
%\ \Gamma_j=\pi(C\circ N^*\s_j),j=1,2.
%\ee
It is easy to verify that $\Gamma_j$ is tangent to $\Gamma_{x_0}$ at $\fp_j$. Then, 
%\be\label{gamma sing corner}
%N^*(\Gamma\cap M_0)=N^*\big((\Gamma_1\cup\Gamma_2\cup\Gamma_{x_0})\cap M_0\big)\cup N^*\fp_1\cup N^*\fp_2.
%\ee
\bs\label{corner wf}
&WF(\hat f_k)\subset N^*\Gamma_1\cup N^*\Gamma_2\cup N^*\Gamma_{x_0}\cup N^*\fp_1\cup N^*\fp_2,\ k=1,2,
\es
in $M_0$. By \eqref{g sing tan} and \eqref{sing g corner},
\bs
&g\in I^{-3/2}(M_0;\Gamma_j)\text{ near any } \fp\in \Gamma_j\setminus\fp_j,\ j=1,2,\\ 
&g\in I^{-2}(M_0;\Gamma_{x_0})\text{ near any } \fp\in \Gamma_{x_0}\setminus\{\fp_1,\fp_2\}.
\es
Hence 
\be\label{WFg c}
WF(g)\subset  N^*\Gamma_1\cup N^*\Gamma_2\cup N^*\Gamma_{x_0}\cup N^*\fp_1\cup N^*\fp_2
\ee
in $M_0$. 
%There will be streaks in $\fr$, but, as before, their primary cause is not a nonlinearity, but the fact that $g$ is not in the range of the Radon transform. 

\begin{figure}[h]
{\centerline{
{\epsfig{file={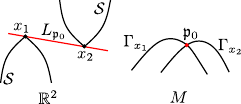}, width=6cm}}
}}
\caption{$L_{\fp_0}$ is tangent to $\s$ in the generalized sense at two corner points.}
\label{fig:two_corner}
\end{figure}

\noindent
{\bf 5.}  Suppose $L_{\fp_0}$ is tangent in the general sense to $\s$ at two points, $x_1$ and $x_2$. The two-corner case is shown in Figure~\ref{fig:two_corner}. As in case 1, $\Gamma\cap M_0=\Gamma_1\cup\Gamma_2$, where $\Gamma_{1,2}$ are smooth and intersect transversely at $\fp_0$. For example, $\Gamma_j=\Gamma_{x_j}$, $j=1,2$, in Figure~\ref{fig:two_corner}. As before, a nonlinear interaction exists at $\fp_0$.  By \eqref{g exp v3}, 
\be\label{WFg gen}
WF(g)\subset N^*\Gamma_1\cup N^*\Gamma_2\cup N^*\fp_0, 
\ee
and, more precisely,  
\be\label{gen dbl case v2}
g\in I^{-r_1,-r_2+(1/2)}(N^*\fp_0,N^*\Gamma_1)+I^{-r_2,-r_1+(1/2)}(N^*\fp_0,N^*\Gamma_2).
\ee
where $r_j=3/2$ in the case of the usual tangency at $x_j$ and $r_j=2$ if $x_j$ is a corner point, $j=1,2$. Both inclusions, \eqref{WFg gen}  and \eqref{gen dbl case v2}, are understood as holding in a neighborhood of $\fp_0$.

\begin{figure}[h]
{\centerline{
{\epsfig{file={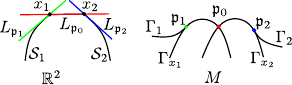}, width=8cm}}
}}
\caption{The line segment case.}
\label{fig:line_segm}
\end{figure}

\noindent
{\bf 6.}  $\s$ contains a flat segment whose endpoints, $x_1$ and $x_2$, are corners, see Figure~\ref{fig:line_segm}. As in case 4, similarly to \eqref{corner wf},
\bs\label{corner-flat wf}
&WF(\hat f_k),WF(g)\subset N^*\Gamma_j\cup N^*\Gamma_{x_j}\cup N^*\fp_j,\ j=1,2,\ k=1,2.
\es
Let $L_{\fp_0}$ be the line that contains the line segment $[x_1,x_2]$.  
Near $\fp_0$, we have similarly to \eqref{gen dbl case v2}
\bs\label{cor cor nl v2}
g\in & I^{-2,-3/2}(N^*\fp_0,N^*\Gamma_{x_1})+I^{-2,-3/2}(N^*\fp_0,N^*\Gamma_{x_2}).
\es

This completes the analysis of possible cases. Now, we prove the remaining claims. 
The above analysis shows that $\Gamma$ is nonsmooth only in case 3 and self-intersects in cases 2, 3, and 6. This proves claims 2 and 3.

From the above analysis, $\Gamma$ is locally smooth only if it corresponds to a smooth segment of $\s$ with a nonvanishing curvature or to a corner point (away from one-sided tangents). Hence, cases 1 and 4 prove claim 4. 

Claim 5 follows from cases 2, 5, and 6 (see \eqref{tang paired lagr v2}, \eqref{gen dbl case v2}, \eqref{cor cor nl v2}).

Claim 6 follows from the results of all six cases, including \eqref{WFg tt}, \eqref{WFg except 1}, \eqref{WFg c}, \eqref{WFg gen}, and \eqref{corner-flat wf}. This completes the proof of Theorem~\ref{thm:main}.

\section{Strength of added singularities.}\label{sec:pp}

There are no added singularities when $L_{\fp_0}$ is tangent to $\s$ at only one point where the curvature of $\s$ is nonzero (case 1 in section~\ref{sec:proof}). Similarly, only the last sum in \eqref{g exp v2} contributes to nonlocal artifacts in the case of tangency at two points (case 2 in section~\ref{sec:proof}). Because we use a linear inversion formula, we separately study the contribution of each term. Pick one of them and denote it
\be
g_a(\al,p)=(p-P_1(\al))_{\iota_1}^{s_1} (p-P_2(\al))_{\iota_2}^{s_2} B(\al,p).
\ee
The subscript ‘$a$’ stands for ‘artifact.’ Because $R^{-1}$ is linear, using a partition of unity, we can assume without loss of generality that $B\in\coi(M_0)$. The reconstruction from $g_a$ uses \eqref{rad inv orig} as follows:
\be\label{rad inv}
f_a(x):=(R^{-1}g_a)(x)=-\frac1{2\pi^2}\int_0^{\pi}\int_{\br} \frac{\pa_p g_a(\al,p)}{p-\vec\al\cdot x}\dd p\dd\al.
\ee

The values $s_1=s_2=1/2$ produce the strongest singularities; therefore, we consider only this case. The other values of $s_1,s_2$ can be analyzed in an analogous manner. 

We begin by setting $\iota_1=\iota_2=+$. The case $\iota_1=\iota_2=-$ is completely analogous. 

\subsection{The $++$ case}
Using entry 26 in the Table of Fourier Transforms of \cite{gs}
\be\label{p+}
p_+^{1/2}=\frac{\Gamma(3/2)e(-3/2)}{2\pi}\int_{\br}(\la-i0)^{-3/2}e^{i\la p}\dd\la,\ e(r)=\exp(ir\pi/2),
\ee 
the 2D Fourier transform of $g_a$ is 
\bs
\tilde g_a(\eta)=\frac{\Gamma^2(3/2)e(-3)}{(2\pi)^2}\int_{M_0}&\int_{\br}\int_{\br}(\la-i0)^{-3/2}(\mu-i0)^{-3/2}B(\al,p)\\ 
&\times e^{i\big[\la(p-P_1(\al))+\mu(p-P_2(\al))+\eta\cdot(\al,p)\big]}\dd\la\dd\mu\dd p\dd\al,
\es
where the integral is understood in the sense of distributions. To compute the Fourier transform, $\tilde g_a(\eta)$, we view $g_a(\al,p)$ as a function on $\br^2$, which is extended by zero outside $M_0$. Set $\eta=r\vec\Theta$, $\vec\Theta=(\Theta_\al,\Theta_p)\in S^1$, $r>0$. Change the variables $\chl=\la/r$ and $\chm=\mu/r$:
\bs
\tilde g_a(r\vec\Theta)=r^2\frac{\Gamma^2(3/2)e(-3)}{(2\pi)^2}&\int_{M_0}\int_{\br}\int_{\br}(r\chl-i0)^{-3/2}(r\chm-i0)^{-3/2}B(\al,p)\\ 
&\times e^{ir\big[\chl(p-P_1(\al))+\chm(p-P_2(\al))+\vec\Theta\cdot(\al,p)\big]}\dd\chl\dd\chm\dd p\dd\al.
\es
The stationary point of the phase is determined by solving the following system:
\bs
&p-P_1(\al)=0,\ p-P_2(\al)=0,\\ 
&-\chl P_1^\prime(\al)-\chm P_2^\prime(\al)+\Theta_\al=0,\ \chl+\chm+\Theta_p=0.
\es
As is easily seen, the stationary point $(\al,p,\chl,\chm)$ is given by
\be\label{stat pt}
\al_0,\, p_0,\, \chl_0(\eta)=\frac{\eta\cdot(1,P_2^\prime(\al_0))}{P_1^\prime(\al_0)-P_2^\prime(\al_0)},\, \chm_0(\eta)=\frac{\eta\cdot(1,P_1^\prime(\al_0))}{P_2^\prime(\al_0)-P_1^\prime(\al_0)},\eta=\vec\Theta.
\ee
The Hessian matrix is given by
\be\label{hessian}
H=\bma 0 & V \\ V^\top & W \ema,\ 
V=\bma -P_1^\prime(\al_0) & 1 \\ -P_2^\prime(\al_0) & 1 \ema,\ 
W=\bma -\chl P_1^{\prime\prime}(\al_0)-\chm P_2^{\prime\prime}(\al_0) & 0 \\ 0 & 0 \ema. 
\ee
The determinant of the Hessian at the stationary point is $(P_1^\prime(\al_0)-P_2^\prime(\al_0))^2\not=0$. Using Sylvester's law of inertia, it is not hard to show that the number of positive and negative eigenvalues of $H$ are the same. By the stationary point method,
\bs\label{g v1}
&\tilde g_a(r\vec\Theta)=C \frac{Q(\vec\Theta)}{r^3}e^{ir\vec\Theta\cdot(\al_0,p_0)}+O(r^{-4}),\ r\to+\infty,\\
%\begin{cases} Q(\vec\Theta)
%,& t\to+\infty,\\
%\overline{Q}(\vec\Theta),& t\to -\infty,\end{cases}\\
&Q(\eta):= e(-3)\big(\chl_0(\eta)-i0\big)^{-3/2}\big(\chm_0(\eta)-i0\big)^{-3/2}
=\kappa|\chl_0(\eta)\chm_0(\eta)|^{-3/2},\\
&C:=\frac{\Gamma^2(3/2)B(\al_0,p_0)}{|P_1^\prime(\al_0)-P_2^\prime(\al_0)|},\ 
\kappa=e\big(-(3/2)[\text{sgn}(\chl_0(\eta))+\text{sgn}(\chm_0(\eta))]\big).
\es
It is easy to verify that $Q(-\eta)={\overline Q}(\eta)$, where the bar denotes the complex conjugation. Also, $\chl_0(\eta)\not=0$ and $\chm_0(\eta)\not=0$ if 
\be\label{eta cond}
\eta\cdot(1,P_j^\prime(\al_0))\not=0,\ j=1,2.
\ee
This assumption is made in what follows. By considering different values of $\text{sgn}(\la_0)$ and $\text{sgn}(\mu_0)$, we compute the corresponding values 
of $\kappa$ using \eqref{g v1}, see Table~\ref{tbl:iota vals}.
\begin{table}[h]
\begin{center}
\begin{tabular}{ c |c| c }
 $\text{sgn}(\la_0)$ & $\text{sgn}(\mu_0)$ & $\kappa$ \\ 
 \hline
 $+$ & $+$ & $i$ \\  
 $+$ & $-$ & $1$ \\   
 $-$ & $+$ & $1$ \\   
 $-$ & $-$ & $-i$ 
\end{tabular}
\end{center}
\caption{The signs of $\chl_0$ and $\chm_0$ and the corresponding values of $\kappa$.}
\label{tbl:iota vals}
\end{table}

In a slightly different form, \eqref{g v1} is
\bs\label{tildeg v1}
\tilde g_a(\eta)= C e^{i \eta \cdot (\al_0,p_0)}Q(\eta)+O(|\eta|^{-4}),\ |\eta|\to\infty.
\es

We select a function $\chi\in\coi(\br^2)$ supported in a sufficiently small neighborhood of $(\al_0,-\vec\al_0^\perp\cdot x_0)$, which equals 1 in a smaller neighborhood of that point.

The reconstruction can be written as follows
\bs\label{fax}
f_a(x)\overset{.}{=} & \frac1{(2\pi)^3}\int_{\br}\int_{\br^2}\chi(\al,\eta_\al/\eta_p) |\eta_p|\tilde g_a(\eta)e^{-i\eta\cdot(\al,\vec\al\cdot x)}\dd\eta\dd\al\\
\overset{.}{=} & \frac C{(2\pi)^3}\int_{\br}\int_{\br}\int_{\br} \chi(\al,\eta_\al/\eta_p)|\eta_p|\big[Q(\eta)+O(|\eta|^{-4})\big]\\
&\hspace{2.5cm}\times e^{-i\eta\cdot(\al-\al_0,\vec\al\cdot x-p_0)}\dd\eta_\al\dd\al\dd\eta_p,
\es
where the above integrals are understood as oscillatory. Here, $\overset{.}{=}$ denotes equality up to a function, which is $C^\infty$ in a small neighborhood of $x_0$. This follows from \eqref{crt canon rel}, \eqref{crtstar canon rel}, and assumption $\text{supp}(g_a)\subset M_0$. Since $x_0\not= x_{1,2}$, it follows that \eqref{eta cond} is satisfied whenever $\chi(\al,\eta_\al/\eta_p)\not=0$.

Suppose $\eta_p>0$. Consider the double integral with respect to $\eta_\al$ and $\al$. Changing variables $\eta_\al\to\check\eta_\al=\eta_\al/\eta_p$ we get
\bs
I_+(\eta_p)= & \eta_p^{-1}\int_{\br^2}\chi(\al,\check\eta_\al)\big[Q(\check\eta_\al,1)+O(\eta_p^{-1})\big]e^{-i\eta_p(\check\eta_\al,1)\cdot(\al-\al_0,\vec\al\cdot x-p_0)}\dd\check\eta_\al\dd\al,
\es
where the big-$O$ term is understood as $\eta_p\to+\infty$. Since $\chi$ is compactly supported, the double integral is over a bounded set. The stationary point is
\be
\al=\al_0,\ \check\eta_\al=-\vec\al_0^\perp\cdot x.
\ee
By the stationary phase method,
\bs\label{I+}
I_+(\eta_p)= & {2\pi Q(-\vec\al_0^\perp\cdot x,1)}{\eta_p^{-2}}e^{-i\eta_p (\vec\al_0\cdot x-p_0)}+O(\eta_p^{-3}),\ \eta_p\to+\infty.
\es

Suppose next $\eta_p<0$. Consider the double integral with respect to $\eta_\al$ and $\al$. Changing variables $\eta_\al\to\check\eta_\al=-\eta_\al/\eta_p$ we get
\bs
I_-(\eta_p)=  |\eta_p|^{-1}\int_{\br}\int_{\br} & \chi(\al,-\check\eta_\al)\big[Q(\check\eta_\al,-1)+O(|\eta_p|^{-1})\big]\\
&\times e^{i(-\eta_p)(-\check\eta_\al,1)\cdot(\al-\al_0,\vec\al\cdot x-p_0)}\dd\check\eta_\al\dd\al.
\es
The stationary point is
\be
\al=\al_0,\ \check\eta_\al=\vec\al_0^\perp\cdot x.
\ee
By the stationary phase method,
\bs\label{I-}
I_-(\eta_p)= & {2\pi Q(\vec\al_0^\perp\cdot x,-1)}{\eta_p^{-2}}e^{-i\eta_p (\vec\al_0\cdot x-p_0)}+O(|\eta_p|^{-3}),\ \eta_p\to-\infty.
\es
Denote
\be
Q_0=Q(-\vec\al_0^\perp\cdot x_0,1).
\ee
By \eqref{stat pt},
\bs\label{stat pt 0}
&\chl_0(-\vec\al_0^\perp\cdot x_0,1)=\frac{P_2^\prime(\al_0)-\vec\al_0^\perp\cdot x_0}{P_1^\prime(\al_0)-P_2^\prime(\al_0)},\\
&\chm_0(-\vec\al_0^\perp\cdot x_0,1)=\frac{P_1^\prime(\al_0)-\vec\al_0^\perp\cdot x_0}{P_2^\prime(\al_0)-P_1^\prime(\al_0)}.
\es
By \eqref{fax}, \eqref{I+}, and \eqref{I-},
\bs\label{Rinv}
f_a(x_0+h\vec\al_0)\sim & \frac{C}{(2\pi)^2}\bigg[Q_0\ioi \eta_p^{-2}e^{-i\eta_p h} \dd\eta_p
+{\overline Q}_0\int_{-\infty}^0 \eta_p^{-2}e^{-i\eta_p h} \dd\eta_p\bigg],
\es
where $\sim$ denotes equality up to a $I^{-3}(N^*0)$ function. Using entries 20 and 28 in the Table of Fourier Transforms of \cite{gs}, the expression in brackets can be written as follows
\bs\label{sing v1}
&(\re Q_0) \CF(\eta_p^{-2})(-h)+i(\im Q_0) \CF(\eta_p^{-2}\text{sgn}\eta_p)(-h)\\
&\sim (\re Q_0)|h|+(\im Q_0)2h\ln|h|.
\es
Therefore the leading singularity of $f_a$ is given by
\bs\label{sing v2}
&f_a(x_0+h\vec\al_0)
\sim C_0\big[(\re \kappa)|h|+(\im \kappa)2h\ln|h|\big],\\
&C_0:=\frac{B(\al_0,p_0)}{16\pi}\frac{|P_1^\prime(\al_0)-P_2^\prime(\al_0)|^2}{\big|(\vec\al_0^\perp\cdot x_0-P_1^\prime(\al_0))(\vec\al_0^\perp\cdot x_0-P_2^\prime(\al_0))\big|^{3/2}},
\es

\begin{figure}[h]
{\centerline{
{\epsfig{file={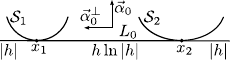}, width=6cm}}
}}
\caption{The case where local segments $\s_1$ (near $x_1$) and $\s_2$ (near $x_2$) of $\s$ are on one side of the double tangent $L_{\fp_0}$.}
\label{fig:tang++}
\end{figure}

\subsection{Example}\label{ssec:pp examp}
To illustrate \eqref{sing v2}, consider the following example. Recall that $x_1$ and $x_2$ are the points at which $L_{\fp_0}$ is tangent to $\s$ (see Figure~\ref{fig:tang++}). Suppose that $\vec\al_0$ points to the half-plane containing the centers of curvature of $\s$ at $x_1$ and $x_2$. Finally, suppose, $\vec\al_0^\perp=(x_1-x_2)/|x_1-x_2|$. We draw $\vec\al_0$ pointing up and $\vec\al_0^\perp$ pointing to the left. Thus, $x_1$ is to the left of $x_2$ and $\vec\al_0^\perp\cdot x_2<\vec\al_0^\perp\cdot x_1$. 

As can be easily seen, $P_2^\prime(\al_0)<P_1^\prime(\al_0)$. This holds because $P_j^\prime(\al_0)=\vec\al_0^\perp\cdot x_j$, $j=1,2$. The last statement is used in case 3 of section~\ref{sec:proof} as the basis for the claim that $\Gamma_j$ is tangent to $\Gamma_{x_0}$ and $\fp_j$, $j=1,2$, see Figure~\ref{fig:corner}. Another variation on the same theme is in case 5 of section~\ref{sec:proof}, see Figure~\ref{fig:line_segm}. Furthermore,
\bs\label{ineqs}
P_2^\prime(\al_0)<P_1^\prime(\al_0)<\vec\al_0^\perp\cdot x_0 &\text{ if }x_0\in (-\infty,x_1),\\
P_2^\prime(\al_0)<\vec\al_0^\perp\cdot x_0<P_1^\prime(\al_0) &\text{ if }x_0\in (x_1,x_2),\\
\vec\al_0^\perp\cdot x_0<P_2^\prime(\al_0)<P_1^\prime(\al_0)  &\text{ if }x_0\in (x_2,\infty).
\es
By \eqref{stat pt 0},
\bs\label{signs}
\text{sgn}(\la_0)=-1,\ \text{sgn}(\mu_0)=+1 &\text{ if }x_0\in (-\infty,x_1),\\
\text{sgn}(\la_0)=-1,\ \text{sgn}(\mu_0)=-1 &\text{ if }x_0\in (x_1,x_2),\\
\text{sgn}(\la_0)=+1,\ \text{sgn}(\mu_0)=-1 &\text{ if }x_0\in (x_2,\infty).
\es

From Table \ref{tbl:iota vals}, the values of $\kappa$ for all $x_0\in L_{\fp_0}\setminus\{x_1,x_2\}$ are given by
\be
\kappa=\begin{cases}
-i,&x_0\in (x_1,x_2),\\
1,& x_0\in L_{\fp_0}\setminus [x_1,x_2].
\end{cases}
\ee
Together with \eqref{sing v2} this implies, up to a positive factor, (see \eqref{B main} and  Figure~\ref{fig:tang++})
\be\label{fa pp final}
f_a(x_0+h\vec\al_0)\sim B_{\frac12, \frac12}(\fp_0)\begin{cases}
h\ln|h|,&x_0\in (x_1,x_2),\\
|h|,& x_0\in L_{\fp_0}\setminus [x_1,x_2].
\end{cases}
\ee

\begin{figure}[h]
{\centerline
{\epsfig{file={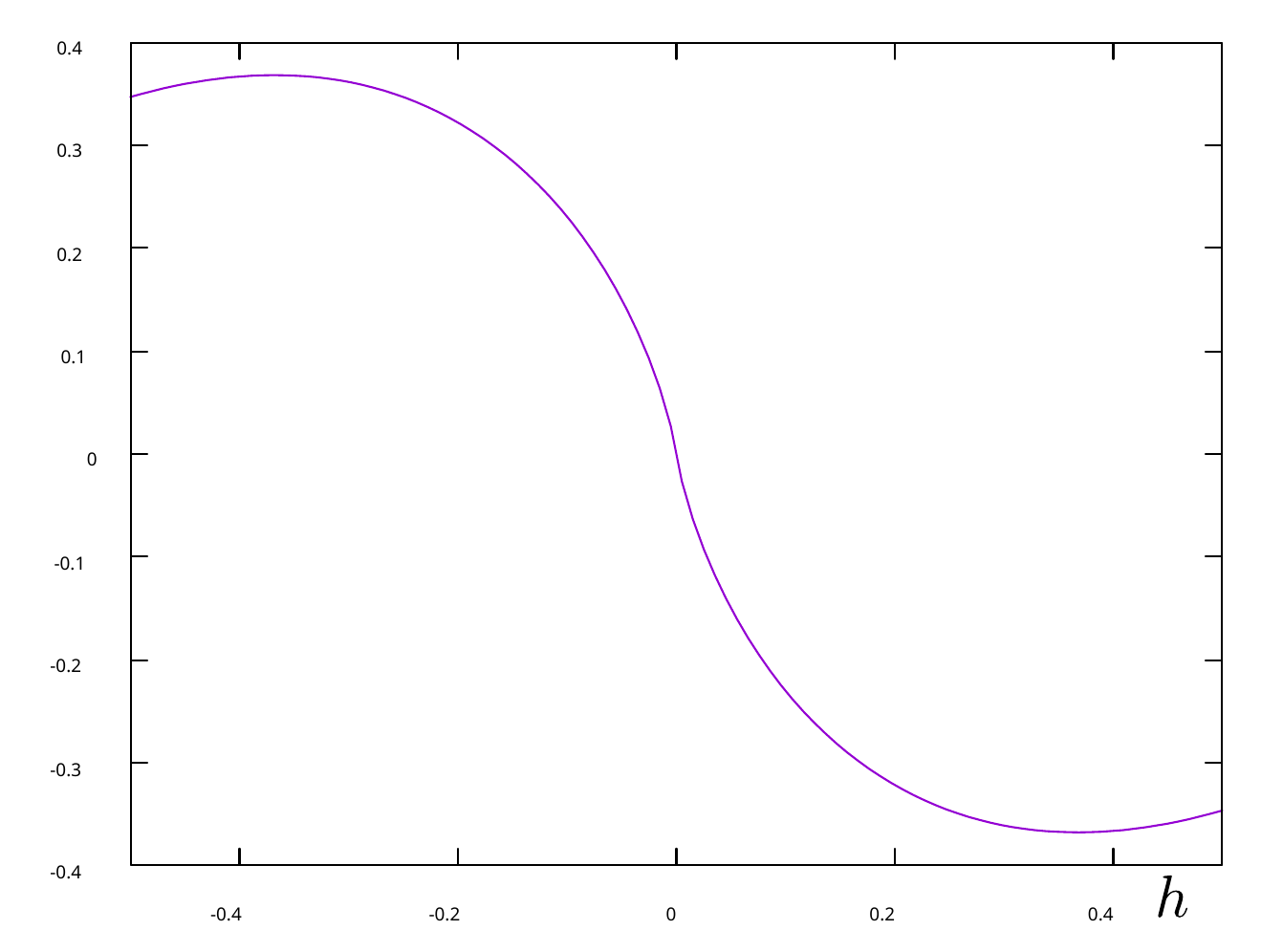}, width=5cm}}}
\caption{Plot of the function $h \log|h|$.}
\label{fig:hlogh}
\end{figure}

For illustration purposes, the graph of the function $h\log|h|$, $|h|\le 0.5$, is shown in Figure~\ref{fig:hlogh}.

%Here we indicated the type 

\section{$+-$ case}\label{sec:pm}

Suppose now that
\be
g_a(\al,p)=(p-P_1(\al))_-^{1/2} (p-P_2(\al))_+^{1/2} B(\al,p),\ B\in\coi(M_0).
\ee
Using \eqref{p+} and entry 25 in the Table of Fourier Transforms of \cite{gs}
\be\label{pminus}
p_-^{1/2}=\frac{\Gamma(3/2)e(3/2)}{2\pi}\int_{\br}(\la+i0)^{-3/2}e^{i\la p}\dd\la,
\ee 
we obtain analogously to \eqref{g v1} and \eqref{tildeg v1}, assuming \eqref{eta cond} holds:
\bs\label{tildeg pm}
\tilde g_a(\eta)= & C Q(\eta)e^{i \eta\cdot(\al_0,p_0)}+O(|\eta|^{-4}),\ |\eta|\to+\infty,\\
Q(\eta):=& \big(\chl_0(\eta)-i0\big)^{-3/2}\big(\chm_0(\eta)+i0\big)^{-3/2}
=\kappa|\chl_0(\eta)\chm_0(\eta)|^{-3/2},\\
\kappa:=&e\big((3/2)[-\text{sgn}(\chl_0(\eta))+\text{sgn}(\chm_0(\eta))]\big),
\es
where $\chl_0(\eta)$ and $\chm_0(\eta)$ are the same as in \eqref{stat pt}, and $C$ is the same as in \eqref{g v1}. As before, $Q(-\eta)={\overline Q}(\eta)$. Therefore, \eqref{tildeg v1}--\eqref{sing v2} still apply, but with the new definition of $Q$. 

\begin{figure}[h]
{\centerline{
{\epsfig{file={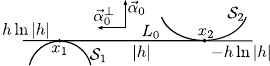}, width=7cm}}
}}
\caption{The case where local segments $\s_1$ (near $x_1$) and $\s_2$ (near $x_2$) of $\s$ are on opposite sides of the double tangent $L_{\fp_0}$.}
\label{fig:tang+-}
\end{figure}

To illustrate this case, we choose $\vec\al_0$ and $\vec\al_0^\perp$ in the same way as before, but $\s_1$ is now below $L_{\fp_0}$ (see Figure~\ref{fig:tang+-}). As is easily seen, \eqref{ineqs} and \eqref{signs} still hold. Substituting \eqref{signs} into the definition of $\kappa$ in \eqref{tildeg pm} we find
\be
\kappa=\begin{cases}
-i,&x_0\in (-\infty,x_1),\\
1,& x_0\in (x_1,x_2),\\
i,&x_0\in (x_2,\infty).
\end{cases}
\ee
Together with \eqref{sing v2} this implies that up to a positive factor one has (see Figure~\ref{fig:tang+-}) 
\be\label{fa pm final}
f_a(x_0+h\vec\al_0)\sim B_{\frac12, \frac12}(\fp_0)\begin{cases}
h\ln|h|,&x_0\in (-\infty,x_1),\\
|h|,& x_0\in (x_1,x_2),\\
-h\ln|h|,&x_0\in (x_2,\infty).
\end{cases}
\ee

\section{Numerical experiments}\label{sec:numerics}

\begin{figure}[h]
{\centerline{\hbox{
{\epsfig{file={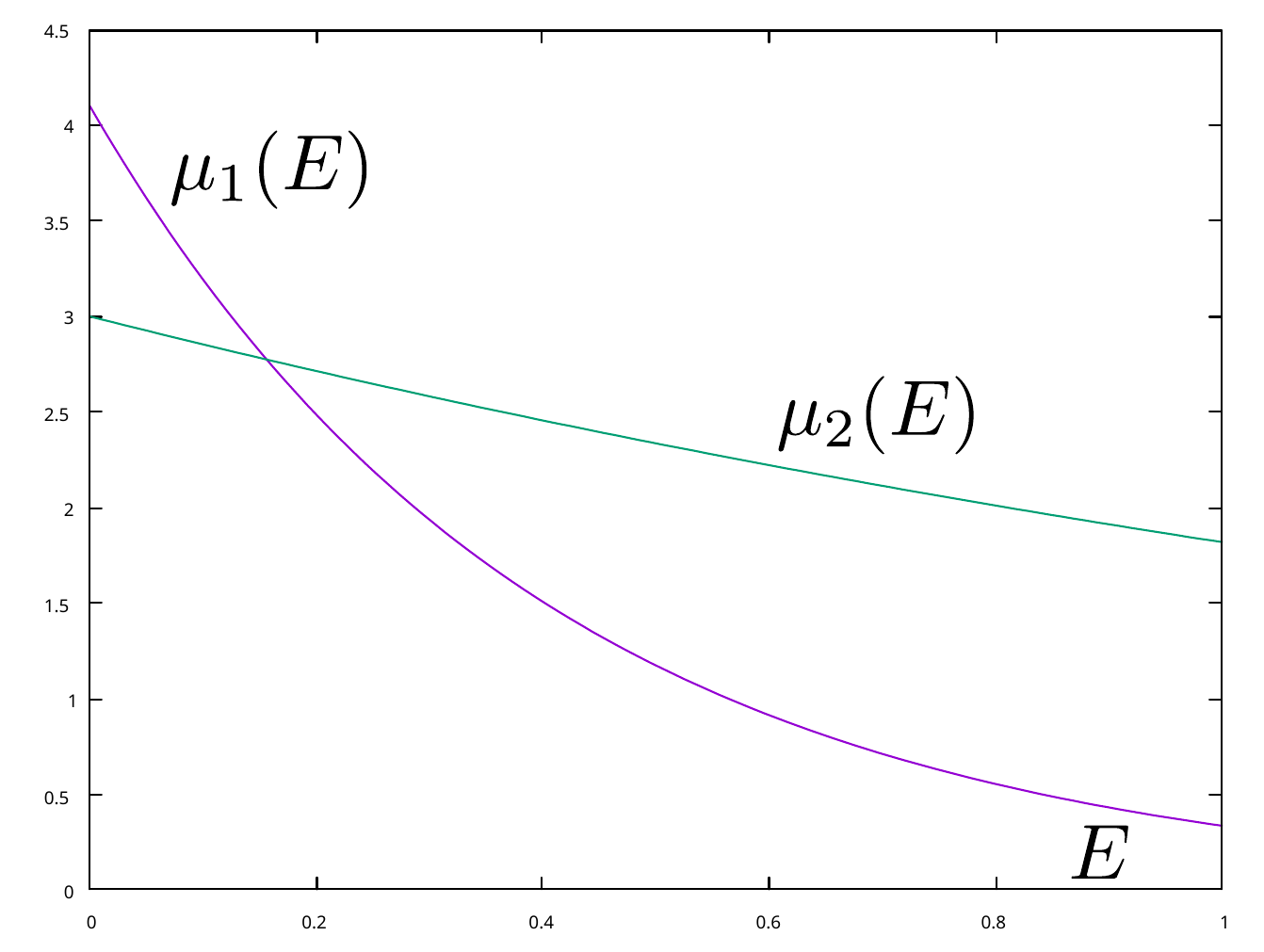}, width=5cm}}
{\epsfig{file={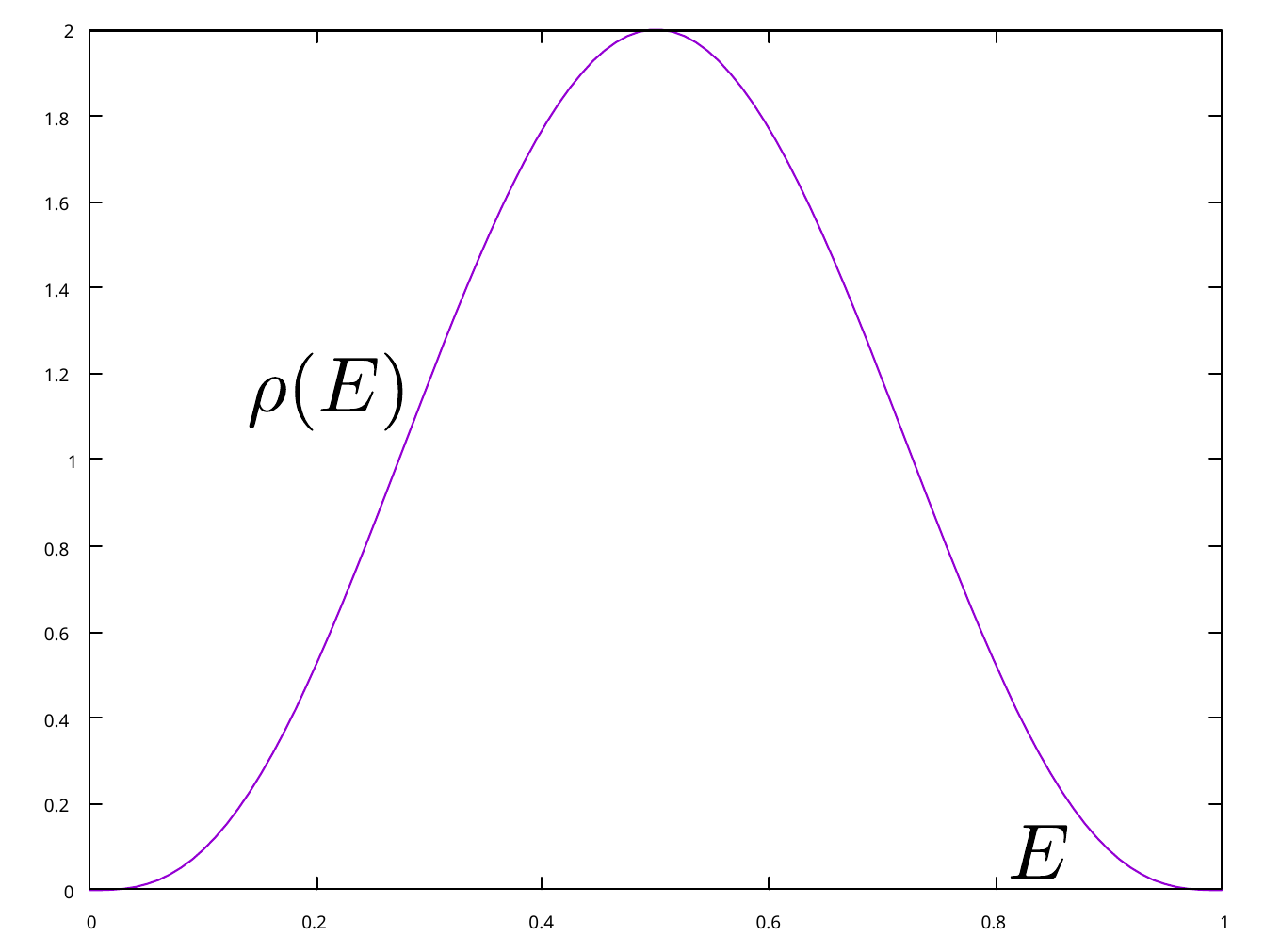}, width=5cm}}
}}}
\caption{Basis functions $\mu_1(E)$ and $\mu_2(E)$ and the energy spectrum $\rho(E)$.}
\label{fig:spec}
\end{figure}

We use a two-ball phantom for the numerical experiments. The two balls have the following centers, $c_k$, and radii, $R_k$:
\be
c_1=(-1.8, -2),\ R_1=1.8,\ c_2=(2.1, 2.2),\ R_2=1.7.
\ee
The two basis functions, $\mu_1$, $\mu_2$, and the normalized product of the X-ray beam energy spectrum and detector efficiency, $\rho(E)$, are given by (see Figure~\ref{fig:spec})
\bs
&\mu_1(E)=4.1\exp(-2.5E),\ \mu_2(E)=3\exp(-0.5E),\\
&\rho(E)=\frac{35}{16}(1-(2E-1)^2)^3,\ 0\le E\le E_{\text{mx}}=1.
\es

Data $G(\al,p)$ are computed using \eqref{data}, where $f_k$ is the characteristic function of the $k$-th ball. We use 
\bs
&N_\al=10^4,\ \Delta_\al=\pi/N_\al,\ N_p=15001,\ \Delta_p=2R_{\text{FOV}}/(N_p-1),\\ 
&\al_k=\Delta_\al k,\ 0\le k<N_\al,\ p_j=-R_{\text{FOV}}+\Delta_p j,\ 0\le j<N_p.
\es
Here $R_{\text{FOV}}$ is the radius of the field of view (FOV).

The reconstruction is computed using \eqref{rad inv}. We denote the reconstructed function $\fr$ instead of $f_a$ as in \eqref{rad inv}, because in computations we use the complete function $g(\fp)$ rather than its paired Lagrangian part $g_a$. The ramp kernel is smoothed by convolving it with a small radius mollifier 
\be
w(t)=\frac{15}{16\e}\big(1-(t/\e)^2\big)^2,\ \e=0.015, 
\ee
to avoid aliasing. 

\begin{figure}[h]
{\centerline{\hbox{
{\epsfig{file={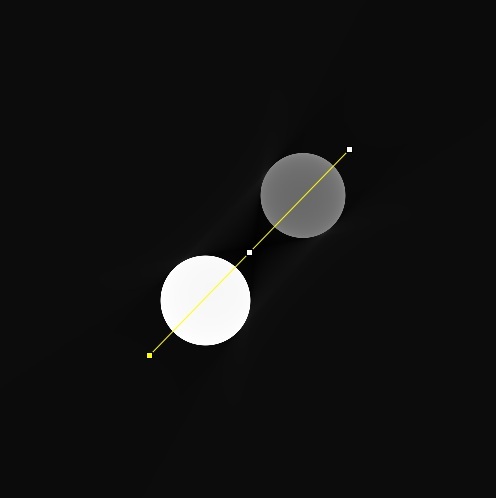}, width=4cm}}
{\epsfig{file={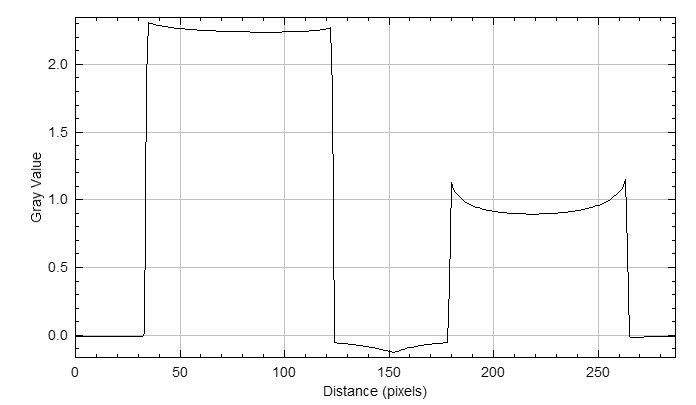}, width=7cm}}
}}}
\caption{Global reconstruction, wide display window width. Typical cupping caused by beam hardening is clearly visible. The plot on the right is along the yellow line shown on the left.}
\label{fig:global wide}
\end{figure}

\begin{figure}[h]
{\centerline{\hbox{
{\epsfig{file={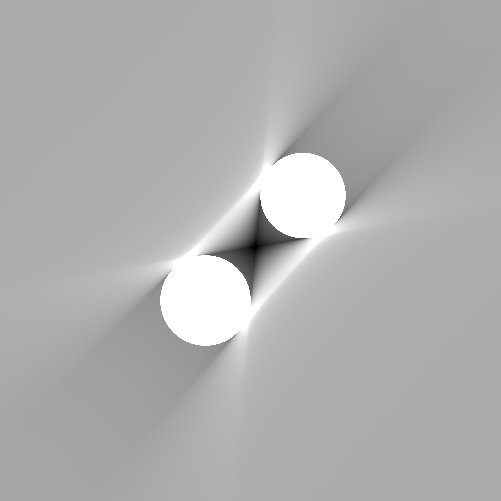}, width=5cm}}
}}}
\caption{Global reconstruction, narrow display window to better see artifacts. WL=-0.03, WW=0.18.}
\label{fig:global narrow}
\end{figure}

Global plots of $\fr(x)$ are shown in Figure~\ref{fig:global wide} and Figure~\ref{fig:global narrow}. Figure~\ref{fig:global wide} shows $\fr(x)$ in the full grayscale window. The right panel in the figure clearly demonstrates low-frequency cupping artifacts, which are caused by beam hardening. The plot on the right is along the yellow line on the left. Figure~\ref{fig:global narrow} shows $\fr(x)$ in the narrow grayscale window $[-0.12,0.06]$. In other words, for display purposes, if $\fr(x)>0.06$ at some $x$, we set $\fr(x)=0.06$. Similarly, if $\fr(x)<-0.12$ at some $x$, we set $\fr(x)=-0.12$. This window is selected to better illustrate the high-frequency BHA. 

\begin{figure}[h]
{\centerline{
{\epsfig{file={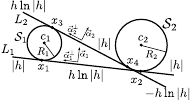}, width=6cm}}
}}
\caption{Illustration of the two-ball phantom used in the numerical experiment. The line $L_1$ corresponds to the $++$ case (section \ref{sec:pp}), and the line $L_2$ - to the $+-$ case (section \ref{sec:pm}).}
\label{fig:diagr}
\end{figure}

Reconstructions are also computed on a dense grid in a neighborhood of several points $x_0$ on the double tangent lines $L_1,L_2$, see Figure~\ref{fig:diagr}. Points $x_0$ on these lines are parameterized by $\nu\in\br$ as follows:
\be
x_0=x_1+\nu(x_2-x_1)\in L_1,\ \nu\not=0,1,
\ee
and similarly for $x_0\in L_2$ using $x_3,x_4$ in place of $x_1,x_2$. To help visualize what our predictions imply for the two-disk phantom we combine Figure~\ref{fig:tang++} and Figure~\ref{fig:tang+-} into a single Figure~\ref{fig:diagr}.

\begin{figure}[h]
{\centerline{\hbox{
{\epsfig{file={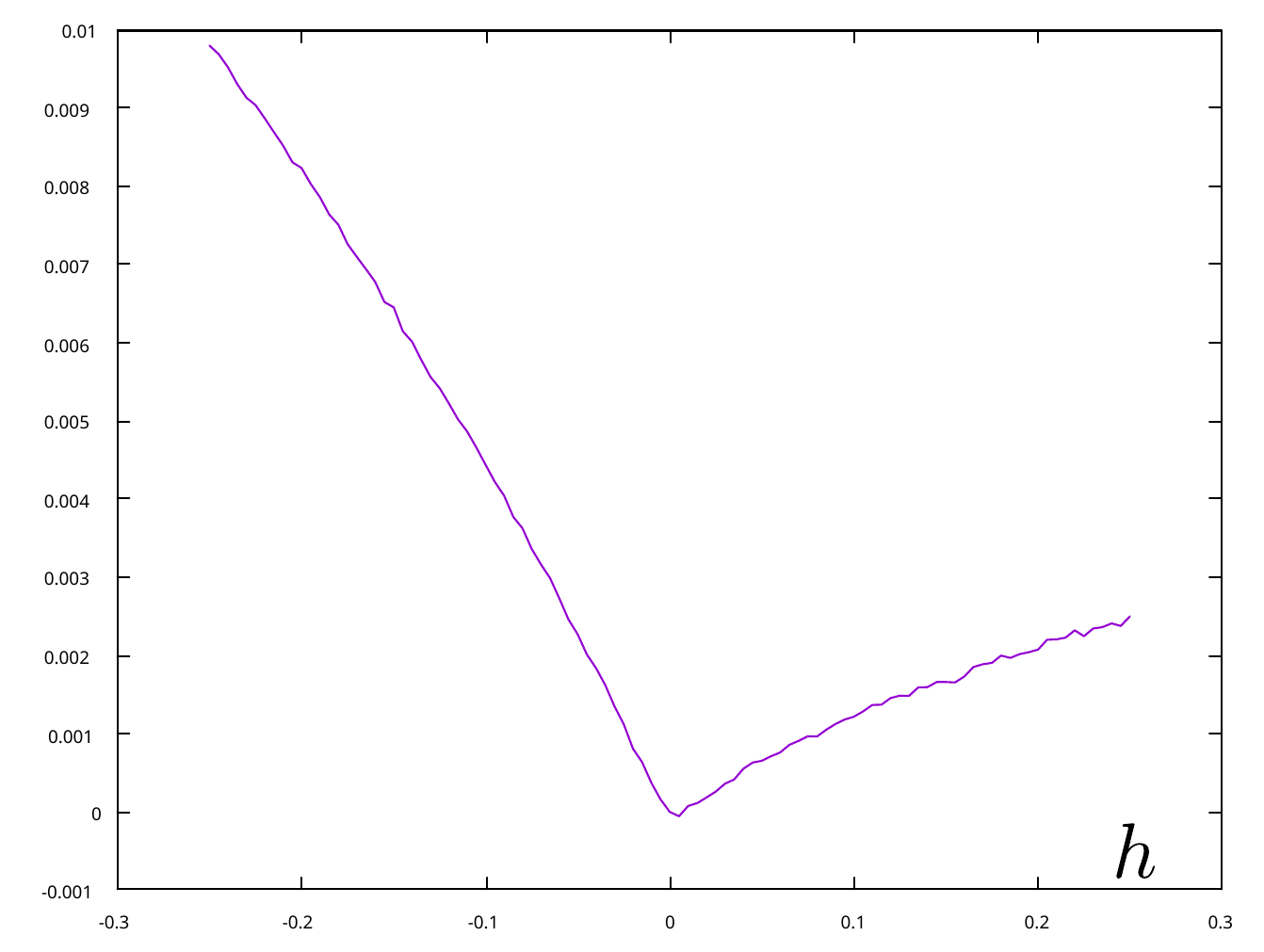}, width=4cm}}
{\epsfig{file={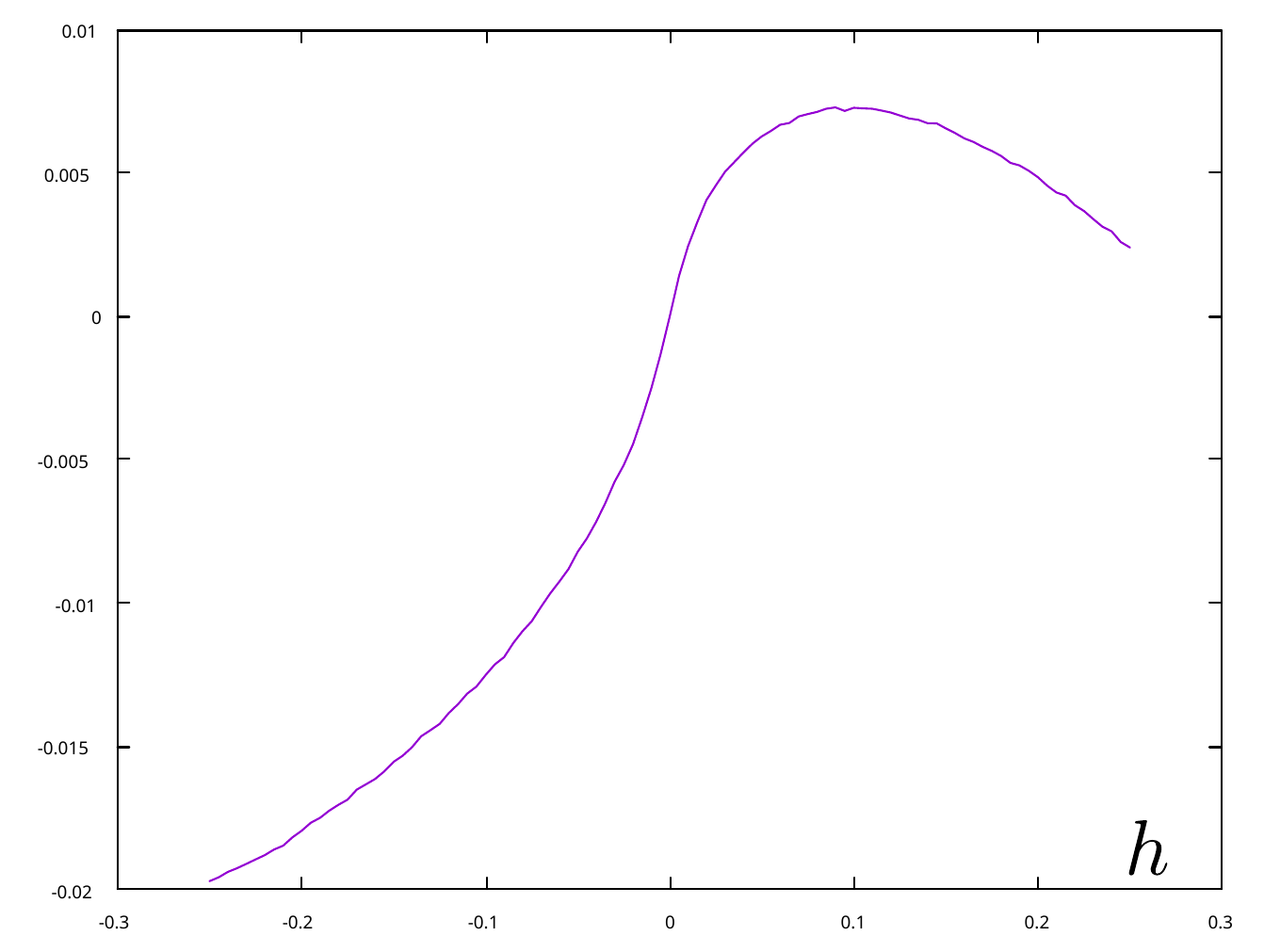}, width=4cm}}
{\epsfig{file={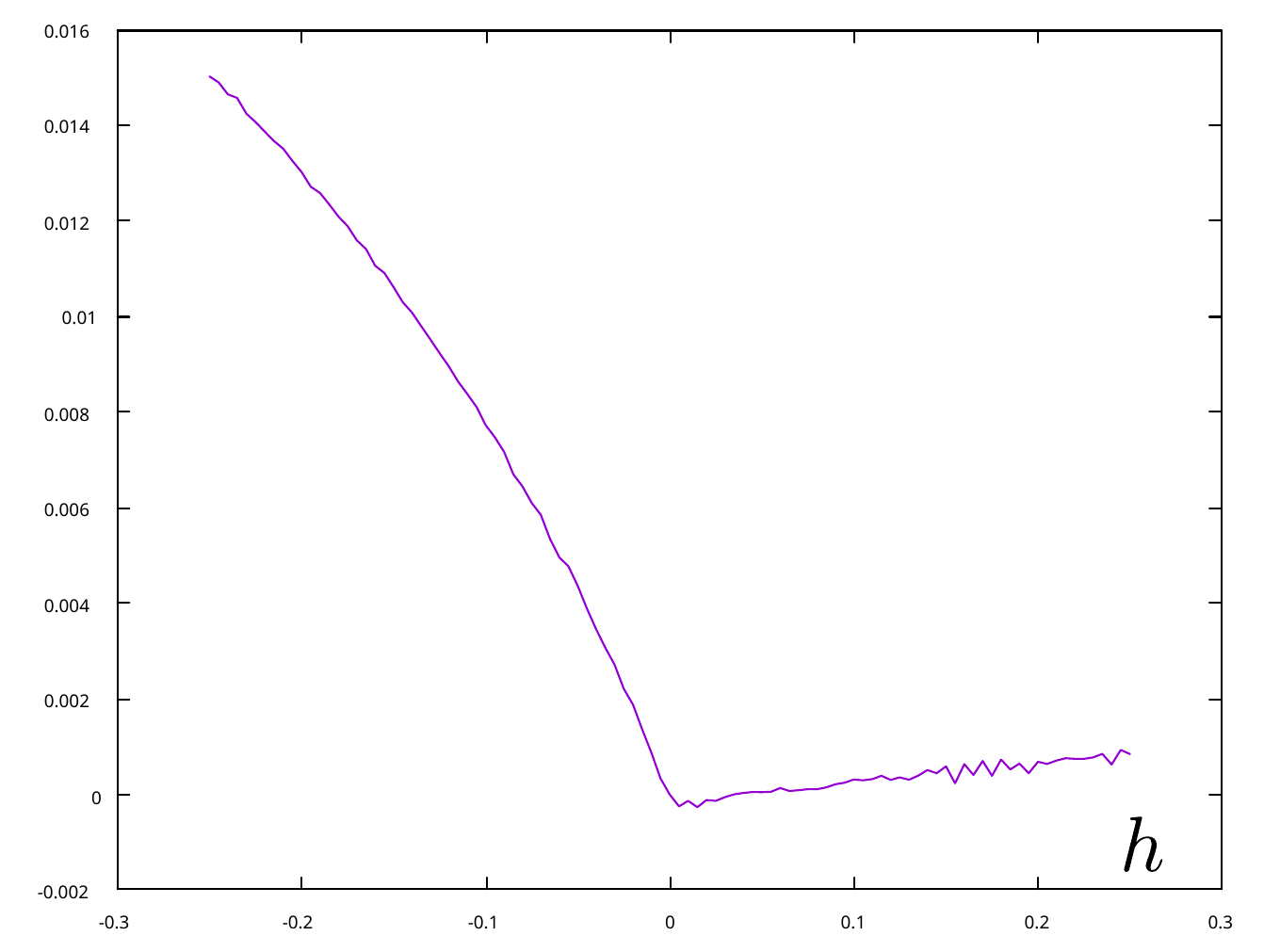}, width=4cm}}
}}}
\caption{Profiles through the artifact at three points on $L_1$, $++$ case. Left: $x_0\in(-\infty,x_1),\nu=-0.4$; middle: $x_0\in(x_1,x_2),\nu=0.4$; right: $x_0\in(x_2,\infty),\nu=1.3$.}
\label{fig:pp art}
\end{figure}

The results of the reconstructions in the $++$ case are shown in Figure~\ref{fig:pp art}. The three plots represent the graphs of $\fr(x)$ across three short-line segments perpendicular to $L_1$ and crossing $L_1$ at three locations corresponding to $\nu=-0.4,0.4,1.3$. The line segments are parameterized as $x=x_0+h\vec\al_1$, $|h|\le 0.25$. Qualitatively, the behavior of $\fr$ matches our  prediction in \eqref{fa pp final}, cf. Figure~\ref{fig:diagr}. The left and right graphs in Figure~\ref{fig:pp art} are shaped like $|h|$, and the middle graph is shaped like $h\ln|h|$ (cf. Figure~\ref{fig:hlogh}). Note that the $|h|$-shaped graphs in the left and right panels point in the same direction, which also agrees with \eqref{fa pp final}. 

In this and the next case, the $|h|$-shaped curves look like a rotated graph of the function $y=|x|$ because of the contribution of smooth terms to $\fr$ that are ignored by our analysis.

\begin{figure}[h]
{\centerline{\hbox{
{\epsfig{file={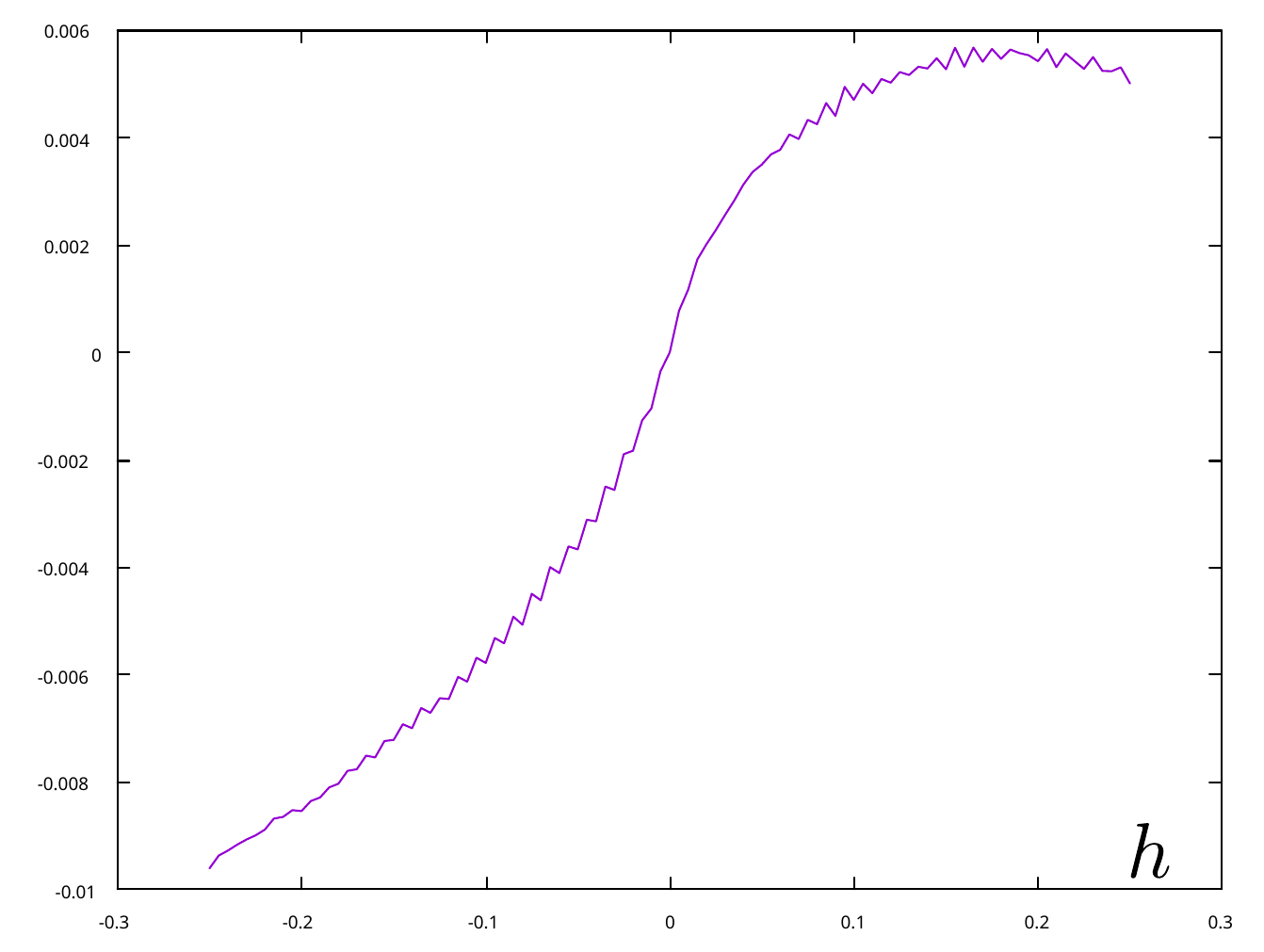}, width=4cm}}
{\epsfig{file={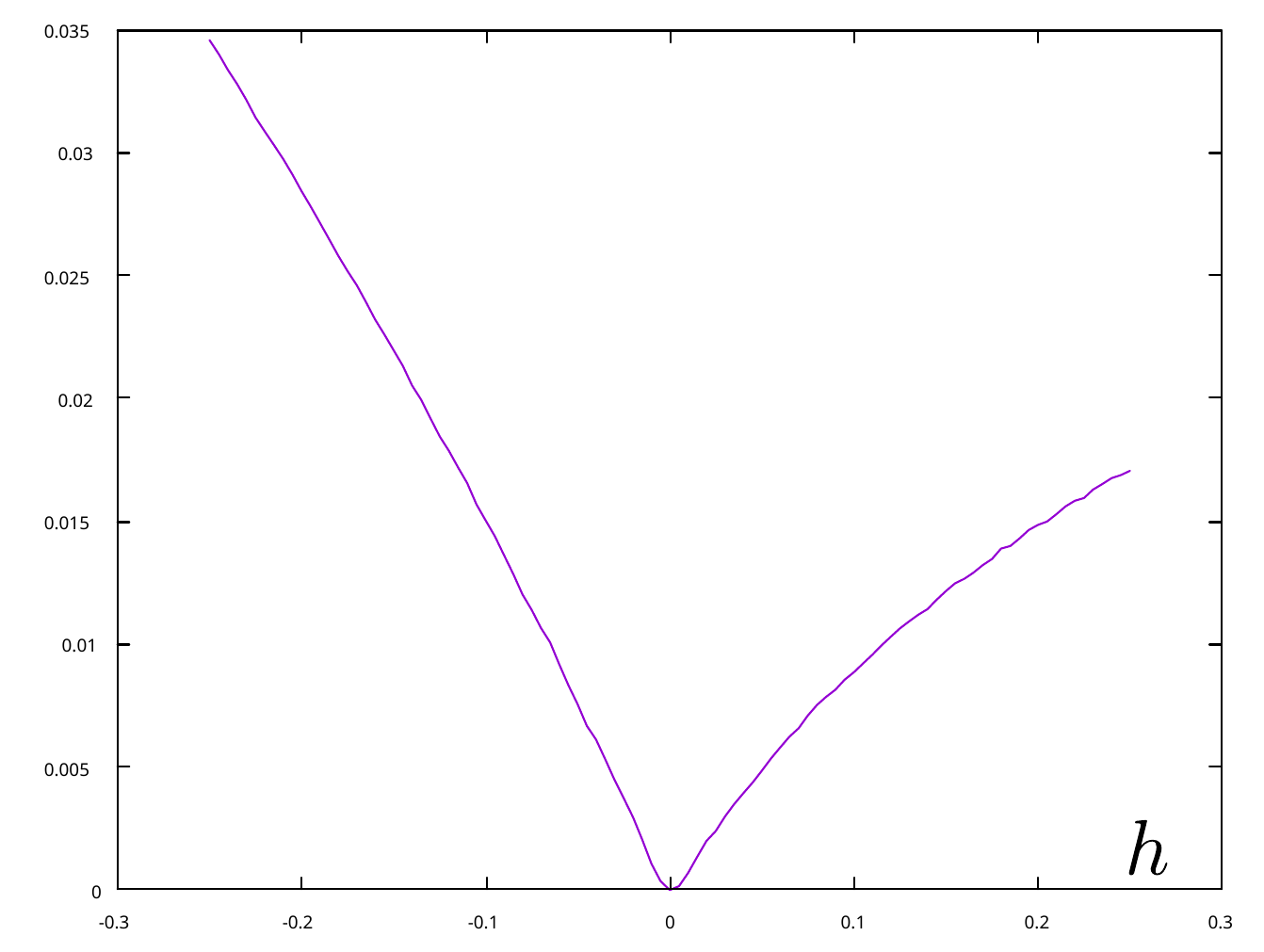}, width=4cm}}
{\epsfig{file={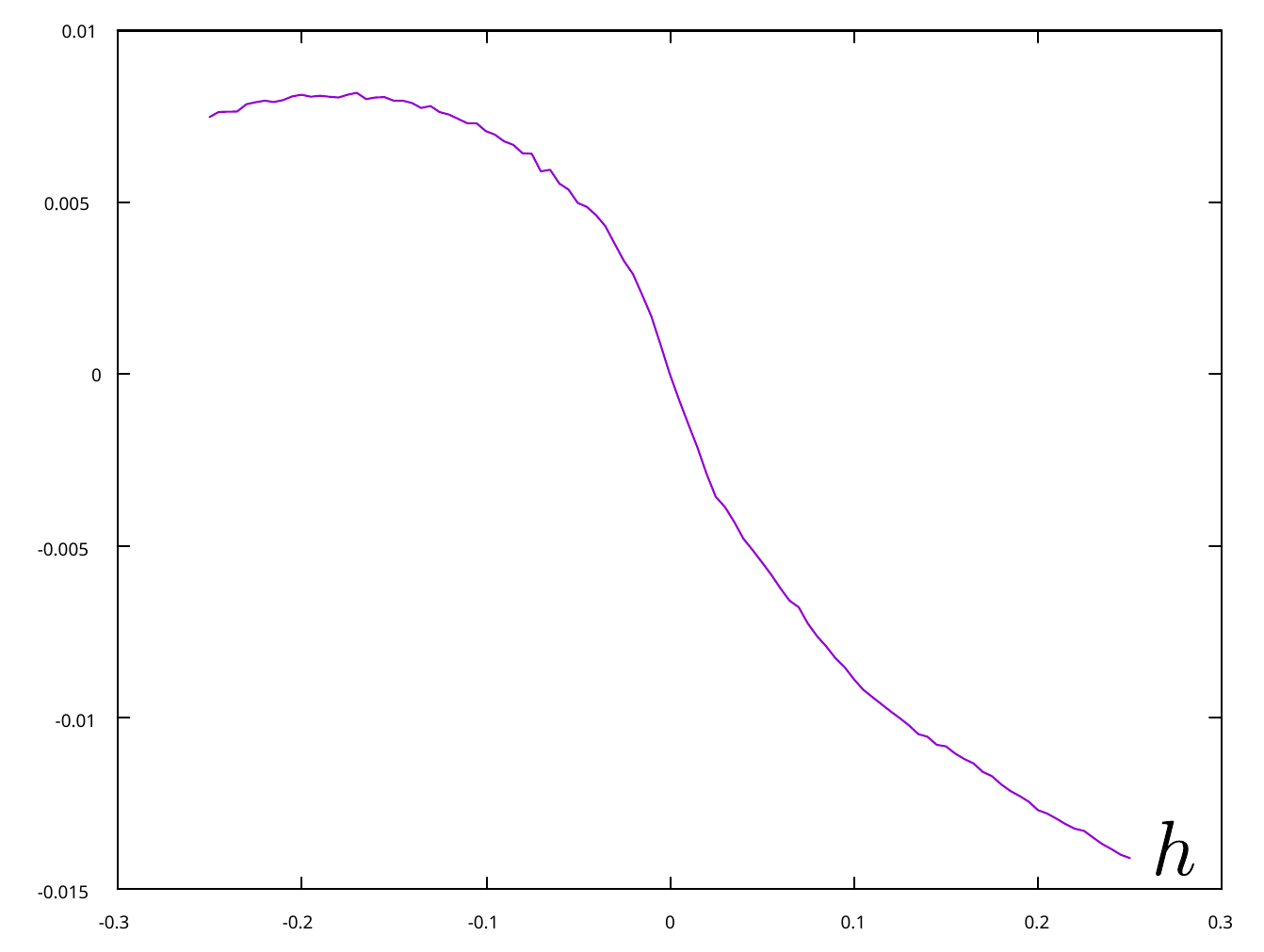}, width=4cm}}
}}}
\caption{Profiles through the artifact at three points on $L_2$, $+-$ case. Left: $x_0\in(-\infty,x_3),\nu=-0.4$; middle: $x_0\in(x_3,x_4),\nu=0.4$; right: $x_0\in(x_4,\infty),\nu=1.3$.}
\label{fig:pm art}
\end{figure}

The results of the reconstructions in the $+-$ case are shown in Figure~\ref{fig:pm art}. The three plots represent the graphs of $\fr(x)$ across three short-line segments perpendicular to $L_2$ and crossing $L_2$ at three locations corresponding to $\nu=-0.4,0.4,1.3$. The line segments are parameterized as $x=x_0+h\vec\al_2$, $|h|\le 0.25$. Qualitatively, the behavior of $\fr$ matches our prediction in \eqref{fa pm final}, cf. Figure~\ref{fig:diagr}. The left and right graphs in Figure~\ref{fig:pm art} are shaped like $h\ln|h|$, and the middle graph is shaped like $|h|$. Note that the $h\ln|h|$-shaped graphs in the left and right panels have prefactors of opposite signs, which is also in agreement with \eqref{fa pm final}.

\bibliographystyle{abbrv}
\bibliography{My_Collection}
\end{document}